\PassOptionsToPackage{svgnames}{xcolor}

\documentclass[preprint]{elsarticle}

\usepackage[margin=0.7in]{geometry}

\usepackage{lineno}

\modulolinenumbers[1]

\usepackage{lipsum}

\usepackage{amsfonts}
\usepackage{graphicx}
\usepackage{epstopdf}
\usepackage{algorithmic}
\usepackage{amssymb}
\usepackage{amsmath}
\usepackage{mathdots} 
\usepackage{graphicx}
\usepackage{lipsum}
\usepackage{amsfonts}
\usepackage{epstopdf}
\usepackage{algorithmic}
\usepackage{algorithm}
\ifpdf
  \DeclareGraphicsExtensions{.eps,.pdf,.png,.jpg}
\else
  \DeclareGraphicsExtensions{.eps}
\fi

\allowdisplaybreaks

\usepackage{todonotes}
\usepackage{multirow}
\usepackage{xcolor}

\usepackage{amsopn}

\newcommand{\bg}{{\mbox{\boldmath{$g$}}}}

\newcommand{\bzeta}{{\mbox{\boldmath{$\zeta$}}}}

\newtheorem{Theorem}{Theorem}

\newtheorem{Definition}{Definition}
\newtheorem{Remark}{Remark}

\newproof{pf}{Proof}
\newcommand{\IRTools}{\textit{IR Tools}}

\usepackage{amsopn}

\begin{document}

\begin{frontmatter}

\title{Flipped structured matrix-sequences in image deblurring with general boundary conditions}

\author[1,2]{Paola Ferrari}

\author[1]{Isabella Furci*}
\cortext[mycorrespondingauthor]{Corresponding author}
\ead{isabella.furci@dima.unige.it}

\author[3,4]{Stefano Serra-Capizzano}

\address[1]{{{Department of Mathematics}}, {University of Genoa}, {{Genoa}, {Italy}}}
\address[2]{{{School of Mathematics and Natural Sciences}}, {University of Wuppertal}, {{Wuppertal}, {Germany}}}
\address[3]{{Department of Science and high Technology}, {University of Insubria}, {{Como}, {Italy}}}
\address[4]{{{Division of Scientific Computing, Department of Information Technology}}, {Uppsala University}, {{Uppsala}, {Swenden}}}

\begin{abstract}

Motivated by a recent work on a preconditioned MINRES for flipped linear systems in imaging, in this note we extend the scope of that research for including more precise boundary conditions such as reflective and anti-reflective ones. We prove spectral results for the matrix-sequences associated to the original problem, which justify the use of the MINRES in the current setting. The theoretical spectral analysis is supported by a wide variety of numerical experiments, concerning the visualization of the spectra of the original matrices in various ways. We also report numerical tests regarding the convergence speed and regularization features of the associated GMRES and MINRES methods. Conclusions and open problems end the present study.

\end{abstract}

\begin{keyword}
Krylov iterative methods, Ill-posedness and regularization problems, eigenvalue distributions,  imaging and signal processing\\65F10, 65F22, 65F15, 94A08
\end{keyword}

\end{frontmatter}

\section{Introduction}
\label{sec:intro}

In recent years a lot of work has been devoted to the study of flipped real nonsymmetric Toeplitz linear systems. The reason is that the flipping (also in a multilevel or even block multilevel setting) leads to a real symmetric indefinite linear system for which an algorithm like MINRES \cite{MR0383715} is well suited: see \cite{MR3323542} for the first algorithmic proposal and for the general idea, where also the preconditioning was proposed and tested.
Subsequently, the global analysis of the asymptotic inertia and of the distribution and clustering were given in \cite{Hon2019, Ferrari2019}, while further results on preconditioning and on the more involved multilevel case were proposed and discussed in \cite{MazzaPestana2018, Ferrari2021, MazzaPestana2021}. In particular, given the indefinite nature of the resulting flipped matrix-sequences, the use of ad hoc positive definite preconditioning leads to preconditioned matrix-sequences clustered at $\pm 1$ in the eigenvalue sense. 
Furthermore we remark that the resulting preconditioned matrices are similar to indefinite symmetric matrices. In this direction, we remind that much earlier a different idea in the indefinite Hermitian Toeplitz setting for obtaining a clustering at $\pm 1$ for the related preconditioned matrix-sequences was presented and studied in \cite{MR1410714}, in connection with a PCG technique. 

On the other hand, this type of structures appear quite frequently in the discretization of evolution equations and hence such kind of technique was also considered in that context \cite{evol-flipped1, evol-flipped2, evol-flipped3} with special attention to PinT methods \cite{Gander-PinT1, Gander-PinT2}. It is worth noticing that in these more applied papers also new theoretical results were deduced concerning the localization of the spectrum so giving a precise bound on the number of iterations of the (preconditioned) MINRES method (see in particular \cite{evol-flipped3}).

Motivated by a recent work on a preconditioned MINRES for flipped linear systems in imaging \cite{donatelli2022symmetrization}, in this note we extend the scope of that research for including more precise boundary conditions such as reflective and anti-reflective ones (see \cite{MR2045058} and references therein). We stress that the idea is the same considered for (block multilevel) Toeplitz matrices with real coefficients and for the approximation of evolution equations, when constant coefficients, Cartesian domains, and uniform grids are taken into account.  In fact, when a space invariant blurring holds, the related structure is $d$-level Toeplitz with real coefficients, $d=1,2,3$, and with (zero) Dirichlet boundary conditions.

As it is transparent from \cite{MR1718798,MR2045058}, the additional difficulty is that the Toeplitz structure is algebraically destroyed by the use of reflective and anti-reflective boundary conditions, even if we prove that the perturbation is spectrally and norm-wise negligible so that it makes sense to apply flipping and then either MINRES or GMRES.
Here the theoretical tools are taken essentially by those arising in the GLT setting, including that of approximating class of sequences. For a general account on the GLT theory and on the related applications see \cite{GSI,GSII,GLT-Eng}.

More precisely, we prove spectral results for the matrix-sequences associated to the original problem, which justify the use of the MINRES/GMRES in the current setting. In fact, the proof of the quasi-Hermitian nature of the related matrix-sequences makes more stable the GMRES while it is also an explanation of the reason why MINRES can be applied. However, in our numerics the two considered Krylov solvers essentially behave in the same way and hence, after an initial comparison, we report only the experiments with the more popular GMRES.

When preconditioning is considered, we briefly discuss how to extend the theory to the case of the preconditioned matrix-sequences. 

The theoretical spectral analysis is supported by a wide variety of numerical experiments, concerning the visualization of the spectra of the original matrices and regarding the convergence speed and regularization features of the considered Krylov solvers.

The paper is organized as follows. Section \ref{sec:notation} is devoted to set the notation and to introduce our spectral machinery. Section \ref{sec:BCs} deals very concisely with reflective and with anti-reflective boundary conditions. Section \ref{sec:main-I} contains the spectral analysis of the nonpreconditioned matrix-sequences, both for reflective and anti-reflective boundary conditions: a short discussion regarding the case of preconditioned matrix-sequences is provided as well. In Section \ref{sec:num} we present and study numerical experiments, while Section \ref{sec:end} contains conclusions and open problems.

\section{Notation and preliminaries}
\label{sec:notation}

Given $X\in\mathbb C^{n\times n}$ and $1\le p\le\infty,\ \|X\|_p$\index{$normxp4$@$\nnorm{X}_p$} denotes the Schatten $p$-norm of $X$, which is defined as the $p$-norm of the vector $(\sigma_1(X),\ldots,\sigma_n(X))$ formed by the singular values of $X$: the Schatten $1$-norm is also called trace norm, while the Schatten $\infty$-norm is referred to as the spectral norm i.e. the maximal singular value. For more details we refer the reader to the book of Bhatia \cite{Bhatia-book}.

We denote by $Y_n$ the backward identity of size $n$ that is

\begin{equation}
Y_n = \begin{bmatrix}{}
1& & &  \\
 & 1 &  &   \\
& & \ddots  &\\
 &  &  & 1
\end{bmatrix} \in \mathbb{R}^{n \times n},
\end{equation} i.e., $[Y_n]_{j,k}=1$ if and only if $j+k=n+1$ and $[Y_n]_{j,k}=0$ elsewhere.

Given two sequences $\{\zeta_n\}_n$ and $\{\xi_n\}_n$, with $\zeta_n\ge0$ and $\xi_n>0$ for all $n$, the notation $\zeta_n=O(\xi_n)$ means that there exists a constant $C$, independent of $n$, such that $\zeta_n\le C\xi_n$ for all $n$; and the notation $\zeta_n=o(\xi_n)$ means that $\zeta_n/\xi_n\to0$ as $n\to\infty$.

\subsection{Image/signal reconstruction}
\label{subs:blur_syst}

Consider a one-dimensional original signal
$$
\tilde{\mathbf{f}}=\left(\ldots, f_{-m+1}, \ldots, f_{0}, f_{1}, \ldots, f_{n}, f_{n+1}, \ldots, f_{n+m}, \ldots\right)^{T}
$$
and the  normalized blurring function given by
$$
\mathbf{h}=\left(\ldots, 0,0, h_{-m}, h_{-m+1}, \ldots, h_{0}, \ldots, h_{m-1}, h_{m}, 0,0, \ldots\right)^{T}, \qquad \sum_{j=-m}^m h_j=1 .
$$
A  blurred signal is the convolution of $\mathbf{h}$ and $\tilde{\mathbf{f}}$, i.e., the $i$-th entry $g_{i}$ of the blurred signal is given by
$$
g_{i}=\sum_{j=-\infty}^{\infty} h_{i-j} f_{j} .
$$
The deblurring problem is to recover the vector $\mathbf{f}=\left(f_{1}, \ldots, f_{n}\right)^{T}$, given the blurring function $\mathbf{h}$ and a blurred signal $\mathbf{g}=\left(g_{1}, \ldots, g_{n}\right)^{T}$ of finite length.
In mathematical terms, having $\mathbf{h}$ and $\mathbf{g}$, we should solve the following underdetermined linear system
\begin{equation}
\label{eq:syst_blur}
\mathbf{g}=\left(\begin{array}{cccccccccc}
h_m & \cdots & h_0 & \cdots & h_{-m}  \\
    & h_m    &     & h_0    &       & h_{-m} &   &   &  0  \\
    &        & \ddots & \ddots & \ddots & \ddots & \ddots \\
    &        &        & \ddots & \ddots & \ddots & \ddots & \ddots \\
    & 0      &        &        & h_m    &     & h_0    &       & h_{-m} \\
    &        &        &        &     & h_m & \cdots & h_0 & \cdots & h_{-m}
\end{array}
\right)
\left(\begin{array}{c}
f_{-m+1} \\
f_{-m+2} \\
\vdots  \\
f_0   \\
f_1   \\
\vdots \\
f_n    \\
f_{n+1} \\
\vdots \\
f_{n+m-1} \\
f_{n+m}
\end{array}
\right).
\end{equation}
Following \cite{MR1718798, MR2045058}, the latter can be rewritten as the sum of three components
$$
T_{l} \mathbf{f}_{l}+T \mathbf{f}+T_{r} \mathbf{f}_{r}=\mathbf{g},
$$
where
$$
\begin{gathered}
T_{l}=\left(\begin{array}{ccc}
h_{m} & \cdots & h_{1} \\
& \ddots & \ddots \\
& & h_{m} \\
0 & &
\end{array}\right), \quad \mathbf{f}_{l}=\left(\begin{array}{c}
f_{-m+1} \\
f_{-m+2} \\
\vdots \\
f_{-1} \\
f_{0}
\end{array}\right), \\
T=\left(\begin{array}{ccccc}
h_{0} & \cdots & h_{-m} & & 0 \\
\vdots & \ddots & \ddots & \ddots & \\
h_{m} & \ddots & \ddots & \ddots & h_{-m} \\
& \ddots & \ddots & \ddots & \vdots \\
0 & & h_{m} & \cdots & h_{0}
\end{array}\right), \quad \mathbf{f}=\left(\begin{array}{c}
f_{1} \\
f_{2} \\
\vdots \\
f_{n-1} \\
f_{n}
\end{array}\right), \\
T_{r}=\left(\begin{array}{cccc}
\vdots \\
h_{-m} & & & \\
\vdots & \ddots & & \\
h_{-1} & \cdots & h_{-m}
\end{array}\right), \quad \text { and } \quad \mathbf{f}_{r}=\left(\begin{array}{c}
f_{n+1} \\
f_{n+2} \\
\vdots \\
f_{n+m-1} \\
f_{n+m}
\end{array}\right) .
\end{gathered}
$$

The use of different boundary conditions (BCs) i.e. zero Dirichlet BCs, nonzero Dirichlet BCs, periodic BCs, reflective BCs, anti-reflective BCs relies on assumptions on the given signal ($d$D images $d\ge 2$) outside the field of values. For the case of $d=1$, we then impose conditions on the values
\[
f_{-m+1}, f_{-m+2}, \ldots, f_{-2}, f_{-1}, f_{n+1}, f_{n+2}, \ldots, f_{n+m-1}, f_{n+m},
\]
in such a way that the system (\ref{eq:syst_blur}) becomes uniquely solvable. It is worthwhile observing that the physics behind the various choices motivates the different denominations used so far. We anticipate that we concentrate our attention on the BCs which ensure the maximal precision in the deblurring process, when piecewise smooth signals for $d=1$ (or $d$D images with $d\ge 2$) are considered. More specifically the reflective BCs ensure the continuity of the $d$D reconstructed signal for every $d\ge 1$ \cite{MR1718798}, while (zero) Dirichlet and periodic BCs fail to lead to continuity in the generic case \cite{MR1718798,MR2045058}, with the important exception of astronomical images with constant (black) background \cite{Book-imaging-astro,Book-imaging}.
Furthermore, as proven in \cite{MR2045058}, anti-reflective BCs ensure the $C^1$ character of the reconstructed signal i.e. for $d=1$, while for $d\ge 2$ the continuity of the reconstructed $d$D signal and that of its normal derivative are guaranteed. More technical issues will be discussed in Section \ref{sec:BCs}.

We conclude the current preliminary subsection on reconstruction problems with some remarks on the 2D case, that is, the application to image restoration. The formal description of the image deterioration is given by the blurring function ${\bf h}$, which in the imaging setting characterizes the Point Spread Function (PSF). Precisely, the PSF is represented as

\[ 
\left[
\begin{array}{cccccccc}
h_{-m,m} &  & \cdots  &  & h_{0,m} &  & \cdots  & h_{m,m} \\
&  &  &  &  &  &  &  \\
\vdots  &  & \ddots  &  & \vdots  &  &  & \vdots  \\
&  &  & h_{-1,1} & h_{0,1} & h_{1,1} &  &  \\
h_{-r,0} &  & \cdots  & h_{-1,0} & {h_{0,0}} & h_{1,0} & \cdots  & h_{r,0} \\
&  &  & h_{-1,-1} & h_{0,-1} & h_{1,-1} &  &  \\
\vdots  &  &  &  & \vdots  &  & \ddots  & \vdots  \\
h_{-m,-m} &  & \cdots  &  & h_{0,-m} &  & \cdots  & h_{m,-m}%
\end{array}%
\right],
\]
where $h_{0,0}$ represent the central coefficient and $2m+1\leq n$. Given a PSF we can define the function

\[
f_{PSF}(\theta_{1},\theta_{2}) = \overset{m}{\underset{j_{1}=-m}{\sum }}\;\underset{%
j_{2}=-m}{\overset{m}{\sum }}h_{j_{1},\,j_{2}}\mathrm{e}^{\iota (j_{1}\theta_{1}+j_{2}\theta_{2})}.
\]


\subsection{Structures and spectral tools}\label{ssec:sp-tools}

The current section is divided into two parts. First we concentrate our attention on block multilevel Toeplitz structures. Then we present spectral tools taken from the relevant literature.

\subsubsection{Preliminaries on block multilevel Toeplitz matrices}
\label{sssec:prelim}

In this section, we provide some useful background knowledge regarding block multilevel Toeplitz matrices.

We let $L^1([-\pi,\pi]^d,\mathbb{C}^{q \times q})$ be the Banach space of all $q\times q$ matrix valued Lebesgue integrable functions over $[-\pi,\pi]^d$, $d,q\ge 1$, equipped with the following norm
\[
\|f\|_{L^1} = \frac{1}{(2\pi)^d}\int_{[-\pi,\pi]^d} \|f(\boldsymbol{\theta})\|_{1}\,{\rm d} \boldsymbol{\theta} < \infty,
\]
where, as recalled at the beginning of the section, $\|A\|_{1}:=\sum_{j=1}^{N}\sigma_{j}(A)$ denotes the trace norm of $A\in \mathbb{C}^{N \times N}$, $\sigma_{j}(A)$, $j=1,\ldots,N$ being the singular values of $A$.

Now we define the multi-index $\mathbf{n}=(n_1,n_2,\dots,n_d)$ where each $n_j$ is a positive integer, $j=1,\ldots,d$.  When writing the expression ${{\bf n}\to\infty}$ we mean that every component of the vector ${\bf n}$ tends to infinity, i.e., $\min_{1\le j\le d} n_j\to \infty$.
{Furthermore, in the current multilevel context, it is convenient to use the Kronecker tensor product $\otimes$ for matrices, where $A\otimes B$ denotes the block matrix of the form $\left(a_{i,j} B\right)$ with $A= \left(a_{i,j}\right)$. In a functional setting,  writing $f=f_1\otimes f_2$ indicates a basic separable function $f(x,y)=f_1(x)f_2(y)$ where $x$ lies in the domain of $f_1$ and $y$ lies in the domain of $f_2$.}

Let $f:$~$[-\pi,\pi]^d\to \mathbb{C}^{q \times q}$ be a function belonging to $L^1([-\pi,\pi]^d,\mathbb{C}^{q \times q})$, and periodically extended to $\mathbb{R}^d$. We define $T_{(\mathbf{n},q)}[f] $ the block multilevel Toeplitz matrix of dimensions $q N(\mathbf{n})\times qN(\mathbf{n})$, with $N(\mathbf{n})= n_1n_2 \dots n_d$ as follows

\begin{equation*}
T_{(\mathbf{n},q)}[f] =\sum_{|j_1|<n_1}\ldots \sum_{|j_d|<n_d} J_{n_1}^{j_1} \otimes \cdots\otimes J_{n_d}^{j_d} \otimes A_{(\mathbf{j})}, \qquad \mathbf{j}=(j_1,j_2,\dots,j_d)\in \mathbb{Z}^d,
\end{equation*}
where
 \[
A_{(\mathbf{j})}=\frac{1}{(2\pi)^d}\int_{[-\pi,\pi]^d}f(\boldsymbol{\theta}){\rm e}^{\iota \left\langle { \bf j},\boldsymbol{\theta}\right\rangle}\, {\rm d}\boldsymbol{\theta},\]
 with $\left\langle { \bf j},\boldsymbol{\theta}\right\rangle=\sum_{t=1}^d j_t\theta_t$, $\iota^2=-1$, are the Fourier coefficients of $f$ and
  $ J^{j}_{n}$  is the $n \times n$ matrix whose $(l,h)$-th entry equals 1 if $(l-h)=j$ and $0$ otherwise.



We indicate by $\{T_{(\mathbf{n},q)}[f]\}_{\mathbf{n}}$ the matrix-sequence whose elements are the matrices $T_{(\mathbf{n},q)}[f]$. The function $f$ is called the \emph{generating function} of $T_{(\mathbf{n},q)}[f]$.


If $f$ is Hermitian-valued almost everywhere (a.e.), then $T_{(\mathbf{n},q)}[f]$ is Hermitian for any choice of $(\mathbf{n}$ and $q$. If $f$ is Hermitian-valued and nonnegative a.e., but not identically zero almost everywhere, then $T_{(\mathbf{n},q)}[f]$ is Hermitian positive definite $\forall \mathbf{n},q$. If $f$ is Hermitian-valued and even a.e., then $T_{(\mathbf{n},q)}[f]$ is real and symmetric for any choice of $\mathbf{n}$ and $q$ \cite{MR2108963,MR2376196,Chan:1996}.
With respect to the previous notation, we have $T_{(\mathbf{n},q)}[f_1\otimes f_2] = T_{(\mathbf{n}[1],q)}[f_1] \otimes T_{(\mathbf{n}[2],q)}[f_2] $ with $\mathbf{n}=(\mathbf{n}[1],\mathbf{n}[2])$, $f_1$ being $k_1$ variate and $f_2$ being $k_2$ variate, with $\mathbf{n}[j]$ being $k_j$ index, $j=1,2$, $k_1+k_2=d$.

\subsubsection{General spectral tools}
\label{sssec:sp-tools}

Throughout this work, we assume that $\phi\in L^1([-\pi,\pi]^d,\mathbb{C}^{q \times q})$ and is periodically extended to $\mathbb{R}^d$. Furthermore, we follow all standard notation and terminology introduced in \cite{GSI} with the only exception of matrix-sizes. In the context of general asymptotic matrix theory, we consider matrix-sequences of the form $\{A_n\}_n$ with $A_n$ of size $d_n$, $d_k<d_{k+1}$, $k\in \mathbb{N}$. Here the matrix-size $d_n$ is always equal to $m N(\mathbf{n})$ as in the block multilevel Toeplitz setting, with $\mathbf{n}=\mathbf{n}(n)$ i.e. all the entries of the vector $\mathbf{n}$ depend on the parameter $n\in \mathbb{N}$ and tend to infinity when $n$ tends to infinity.

We are now ready for introducing our spectral tools specialized to our context.


\begin{Definition}{\rm (Eigenvalue and singular value distribution of a matrix-sequence)}\label{def:distributions}
Let $C_c(\mathbb{C})$ (or $C_c(\mathbb{R})$) be the space of complex-valued continuous functions defined on $\mathbb{C}$ (or $\mathbb{R}$) with bounded support. Let $\{A_{(\mathbf{n},q)}\}_{n}$ be a matrix-sequence with $A_{(\mathbf{n},q)}$ of size $q N(\mathbf{n})$ and $\mathbf{n}=\mathbf{n}(n)$, $n\in \mathbb{N}$.
\begin{enumerate}
\item We say that $\{A_{(\mathbf{n},q)}\}_{{n}}$ has an asymptotic singular value distribution described by a $q\times q$ matrix-valued $\phi$ if 
\begin{multline*}
\lim_{\mathbf{n}\to \infty} \frac{1}{q N(\mathbf{n})}\sum_{j=1}^{q N(\mathbf{n})}F(\sigma_j(A_{(\mathbf{n},q)}))=\\\frac{1}{(2\pi)^d}\int_{[-\pi,\pi]^d} \frac{1}{q} \sum_{j=1}^{q} F( \sigma_j(\phi(\mathbf{x}))  )\,{\rm d}\mathbf{x},\quad \forall F \in C_c(\mathbb{R}).
\end{multline*}
\item We say that $\{A_{(\mathbf{n},q)}\}_{{n}}$ has an asymptotic eigenvalue (or spectral) distribution described by a matrix-valued $\phi$ if 
\begin{multline*}
\lim_{\mathbf{n}\to \infty} \frac{1}{q N(\mathbf{n})}\sum_{j=1}^{q N(\mathbf{n})}F(\lambda_j(A_{(\mathbf{n},q)}))=\\ \frac{1}{(2\pi)^d}\int_{[-\pi,\pi]^d} \frac{1}{q} \sum_{j=1}^{q} F( \lambda_j(\phi(\mathbf{x}))  )\,{\rm d}\mathbf{x},\quad \forall F \in C_c(\mathbb{C}).
\end{multline*}
\end{enumerate}
\end{Definition}

A notion which is strictly related to that of distribution is the concept of clustering: they are both inherently asymptotic notions and the first reduces to the second in special cases. Given $z\in\mathbb C$, $\epsilon>0$, the notation $B(z,\epsilon)$ denotes the closed disk with center $z$ and radius $\epsilon$, $B(z,\epsilon)\doteq\{w\in\mathbb C:\,|w-z|<\epsilon\}$. For $S\subseteq\mathbb C$ and $\epsilon>0$, we the notation $B(S,\epsilon)$ indicates the $\epsilon$-expansion of $S$, that is $B(S,\epsilon)\doteq\bigcup_{z\in S}B(z,\epsilon)$.
\begin{Definition}\label{def-cluster}
  Let $\{X_n\}_n$ be a sequence of matrices, with $X_n$ of size $d_n=q N(\mathbf{n})$  with $\mathbf{n}=\mathbf{n}(n)$, $n\in \mathbb{N}$, and
  let $S\subseteq\mathbb C$ be a nonempty closed subset of $\mathbb C$. $\{X_n\}_n$ is {\em strongly clustered} at $S$
  in the sense of the eigenvalues if, for each $\epsilon>0$, the number of eigenvalues of $X_n$ outside $B(S,\epsilon)$
  is bounded by a constant $q_\epsilon$ independent of $n$. In symbols,
  $$q_\epsilon(n,S)\doteq\#\{j\in\{1,\ldots,d_n\}: \lambda_j(X_n)\notin B(S,\epsilon)\}=O(1),\quad\mbox{as $n\to\infty$.}$$
  $\{X_n\}_n$ is {\em weakly clustered} at $S$ if, for each $\epsilon>0$,
  $$q_\epsilon(n,S)=o(d_n), \quad \mbox{as $n\to\infty.$}$$
  If $\{X_n\}_n$ is strongly or weakly clustered at $S$ and $S$ is not connected, then the connected components of $S$ are called sub-clusters.
\end{Definition}
 We recall that, for a measurable function $g:D\subseteq\mathbb R^d\to\mathbb C$, $D=[-\pi,\pi]^d$, the essential range of $g$ is defined as
  $\mathcal{ER}(g)\doteq\{z\in\mathbb C:\,\mu_t(\{g\in
  B(z,\epsilon)\})>0\mbox{ for all $\epsilon>0$}\}$, where $\{g\in
  B(z,\epsilon)\}\doteq\{x\in D:\,g(x)\in B(z,\epsilon)\}$.
  $\mathcal{ER}(g)$ is always closed and if $g$ is continuous and $D$ is contained in the closure of its interior, then
  $\mathcal{ER}(g)$ coincides with the closure of the image of $g$.   Hence, if $\{X_n\}_n\sim_\lambda \phi$ (with $\{X_n\}_n,\  \phi$ as in Definition \ref{def:distributions}), then, by \cite[Theorem 4.2]{MR2287378}, $\{X_n\}_n$ is weakly clustered at the essential range of $\phi$, defined
  as the union of the essential ranges of the eigenvalue functions $\lambda_i(\phi),\ i=1,\ldots, m$: $\mathcal{ER}(\phi)\doteq\bigcup_{i=1}^m\mathcal{ER}(\lambda_i(\phi))$. Finally, if $\mathcal{ER}(\phi)=s$ with $s$ fixed complex number and $\{X_n\}_n\sim_\lambda \phi$, then $\{X_n\}_n$ is weakly clustered at $s$.
  All the considerations above can be translated in the singular value setting as well, with obvious minimal modifications.
 For further relationships among distribution, clustering, Schatten $p$-norms the reader is referred to \cite{taud2} (see also \cite{MR1366576}).

Finally, we introduce the following definitions and a key lemma regarding the notion of approximating class of sequences, in order to prove our main distribution results in the next {section}.
	
	\begin{Definition}{\rm (Approximating class of sequences)}\label{def:ACS}
		Let $\{A_{(\mathbf{n},q)}\}_{{n}}$ be a matrix-sequence with $A_{(\mathbf{n},q)}$ of size $d_n=q N(\mathbf{n})$, $\mathbf{n}=\mathbf{n}(n)$, $n\in \mathbb{N}$, and let $\{\{B_{{(\mathbf{n},q)},j}\}_{{n}}\}_j$ be a class of matrix-sequences of the same matrix sizes. We say that $\{\{B_{{(\mathbf{n},q)},j}\}_{{n}}\}_j$ is an \textit{approximating class of sequences (a.c.s.)} for $\{A_{(\mathbf{n},q)}\}_{{n}}$ if the following condition is met: for every $j$ there exists $\mathbf{n}_j$ such that, for $  \mathbf{n}\ge\mathbf{n}_j$,
		\[
		A_{(\mathbf{n},q)}=B_{{(\mathbf{n},q)},j}+R_{{(\mathbf{n},q)},j}+N_{{(\mathbf{n},q)},j},
		\]
		\[
		\textrm{rank}~R_{{(\mathbf{n},q)},j}\leq c(j)N(\mathbf{n}) \quad {\rm and} \quad \|N_{\mathbf{n},q}\|\leq\omega(j),
		\]
		where $\mathbf{n}_j$, $c(j)$, and $\omega(j)$ depend only on $j$ and \[\lim_{j\to\infty}c(j)=\lim_{j\to\infty}\omega(j)=0.\]
	\end{Definition}
	
	We use $\{\{B_{{(\mathbf{n},q)},j}\}_{{n}}\}_j \xrightarrow{\text{a.c.s.\ wrt\ $j$}}\{A_{(\mathbf{n},q)}\}_{{n}}$ to denote that $\{\{B_{{(\mathbf{n},q)},j}\}_{\mathbf{n}}\}_j$ is an a.c.s. for $\{A_{(\mathbf{n},q)}\}_{\mathbf{n}}$.
	
	The following is a useful criterion to identify an a.c.s. without constructing the splitting present in Definition \ref{def:ACS}, while the subsequent one is a powerful approximation theorem for identifying the distribution of an involved matrix-sequence starting from simpler matrix-sequences.
	\begin{Theorem}\label{thm:acs_caratt}
	Let $\{A_{(\mathbf{n},q)}\}_{{n}}$ be a matrix-sequence with $A_{(\mathbf{n},q)}$ of size $d_n=q N(\mathbf{n})$, $\mathbf{n}=\mathbf{n}(n)$, $n\in \mathbb{N}$, and let
	$\{\{B_{{(\mathbf{n},q)},j}\}_{{n}}\}_j$ be a class of matrix-sequences of the same matrix sizes. Suppose that for every $j$ there exists $\mathbf{n}_j$ such that, for $  \mathbf{n}\ge\mathbf{n}_j$
	\[\|A_{(\mathbf{n},q)}- B_{{(\mathbf{n},q),j}}\|_{1}\le \epsilon(j)N(\mathbf{n}),\]
where $\lim_{j\to \infty}  \epsilon(j)= 0$. Then,
\[
\{\{B_{{(\mathbf{n},q)},j}\}_{{n}}\}_j\xrightarrow{\text{a.c.s.\ wrt\ $j$}}\{A_{(\mathbf{n},q)}\}_{{n}}.
\]
	\end{Theorem}
\begin{Theorem}\label{lem:Corollary5.1}
		Let $\{A_{(\mathbf{n},q)}\}_{{n}}, \{B_{{(\mathbf{n},q)},j}\}_{{n}}$ $j\in \mathbb{N}$, $\mathbf{n}=\mathbf{n}(n)$, $n\in \mathbb{N}$, be matrix-sequences and let $\phi,\phi_j:D \subset \mathbb{R}^d \to \mathbb{C}$, $D=[-\pi,\pi]^d$, be measurable functions defined on $D$. Suppose that
		
		\begin{enumerate}
			\item $\{\{B_{{(\mathbf{n},q)},j}\}_{{n}}\}_j \sim_{\sigma}  \phi_j$ for every $j$,
			\item $\{\{B_{{(\mathbf{n},q)},j}\}_{{n}}\}_j \xrightarrow{\text{a.c.s.\ wrt\ $j$}} \{A_{(\mathbf{n},q)}\}_{{n}}$,
			\item $\phi_j \to \phi$ in measure.
		\end{enumerate}
		
		Then
		\[
		\{A_{(\mathbf{n},q)}\}_{{n}} \sim_{\sigma}  \phi.
		\]
		
		Moreover, if the first assumption is replaced by $\{B_{{(\mathbf{n},q)},j}\}_{\mathbf{n}}\sim_{\lambda}  \phi_j$ for every $j$, given that the other two assumptions are left unchanged, and all the involved matrices are Hermitian, then  $\{A_{(\mathbf{n},q)}\}_{{n}} \sim_{\lambda}  \phi$.
	\end{Theorem}

Before discussing the asymptotic spectral distribution of the resulting matrix-sequences, it is crucial to our preconditioning theory development that we introduce the following notation:

	\begin{Definition}{\rm \cite[Definition 3.1]{Ferrari2021}}\label{def:psi}
	Given the vector $\mathbf{p}=[2\pi,2\pi,\dots,2\pi]^T\in \mathbb{R}^d$ and a function $g$ defined over $[0,2\pi]^d$, we define
$\psi_{g}$ over $[-2\pi,0]^d\cup [0,2\pi]^d$ in the following manner
	\begin{equation}\label{eq:psi}
\psi_g(\boldsymbol{\theta})=\left\{
	\begin{array}{cc}
	g(\boldsymbol{\theta}), & \boldsymbol{\theta}\in [0,2\pi]^d, \\
	-g(\boldsymbol{\theta}+\mathbf{p}), & \boldsymbol{\theta}\in  [-2\pi,0]^d, \ \boldsymbol{\theta}\neq \mathbf{0}.
	\end{array}
	\right.\,
	\end{equation}

	\end{Definition}
	
	This is a general tool useful for the latter purpose and has a crucial role in the description of the spectrum of the involved matrix-sequences, which was first established in \cite{Ferrari2019} for symmetrized block Toeplitz matrix-sequences and then used in \cite{Ferrari2021} for symmetrized multilevel Toeplitz matrix-sequences. We also refer the detailed analyses on the same problems in the uni-level case \cite{MazzaPestana2018} and the multilevel case \cite{MazzaPestana2021}. It should be observed that in the latter two works the distribution function is defined as
	 \begin{equation}\label{eq:psi-bis}
\psi_g(\boldsymbol{\theta})={\rm diag}\left(g(\boldsymbol{\theta}),-g(\boldsymbol{\theta})\right),\ \ \ \boldsymbol{\theta}\in [-\pi,\pi]^d.
	\end{equation}
This non-uniqueness of the spectral or singular value symbol is well known \cite{ES-NLAA,rearr}. In fact, any rearrangement of a symbol is still a symbol and the two representations in (\ref{eq:psi}) and (\ref{eq:psi-bis}) are equivalent in this respect. However, we point out that the representation in (\ref{eq:psi-bis}) fits better with Definition \ref{def:distributions}, because of our technical assumption of fixing the definition domain of the symbol as $[-\pi,\pi]^d$.

\begin{Theorem}{\cite[Corollary 4.1]{GSI}\cite[Corollory 2.3]{GSII}}
\label{th:GLT-lambda}
Let $\{X_n\}_n,\,\{Z_n\}_n$ be sequences of matrices, with $X_n$, $Z_n$ of size $d_n$,  and set $A_n=X_n+Z_n$, with $d_n=mN(\mathbf{n})$, $\mathbf{n}=\mathbf{n}(n)$, $n\in \mathbb{N}$. Assume that the following conditions are met.
\begin{enumerate}
	\item $\|X_n\|_\infty,\|Z_n\|_\infty\le C$ for all $n$, where $C$ is a constant independent of $n$.
  \item Every $X_n$ is Hermitian and $\{X_n\}_n\sim_\lambda \kappa$.
  \item $\|Z_n\|_1=o(d_n)$.
\end{enumerate}
Then $\{A_n\}_n\sim_\lambda \kappa$.

Furthermore, $\|Z_n\|_1=O(1)$ then the range of $\kappa$ is a strong cluster for the eigenvalues of $\{A_n\}_n$.
\end{Theorem}

\begin{Remark}
Using the tools in \cite{MR2287378}, if the assumptions in Theorem \ref{th:GLT-lambda} are satisfied with $\|Z_n\|_1=O(1)$, then the essential range of $\kappa$, $\mathcal{ER}(\kappa)$, is a strong cluster for the eigenvalues of $\{A_n\}_n$. That is, $\forall \epsilon >0, $ $\exists \, c_\epsilon, \, n_\epsilon$ such that
\[\# \{\lambda_j\left(A_n\right)\, : \,  {\rm dist} \left(\lambda_j\left(A_n\right), \mathcal{ER}(\kappa) \right)>\epsilon\}\le c_\epsilon, \, \forall n\ge n_\epsilon.\]
\end{Remark}
\begin{Remark}\label{going higher}
In the rest of the paper we have $d=1,2$ and $q=1$. However, we presented all the mathematical tools in more generality since the case $d=3$ occurs in $3D$ imaging, while $q\ge 2$ occurs in several imaging contexts such as missing data, super-resolution, inpainting, colors etc (see e.g \cite{qlarger1, qlarger2, qlarger3} and references therein). We refer to the latter in the conclusions. 
\end{Remark}

\section{Reflective and anti-reflective boundary conditions}
\label{sec:BCs}
In Subsection \ref{subs:blur_syst} we mentioned that a blurred signal $\mathbf{g}$ is determined not by $\mathbf{f}$ only, but also by
$$
\left(f_{-m+1}, \ldots, f_{0}\right)^{T} \quad \mbox{ and } \quad \left(f_{n+1}, \ldots, f_{n+m}\right)^{T}
$$
and consequently the linear system (\ref{eq:syst_blur}) is underdetermined.  Then,
when modelling the blur of a signal, it is crucial to make certain assumptions (called boundary conditions) on the unknown data $f_{-m+1}, \ldots, f_{0}$ and $f_{n+1}, \ldots, f_{n+m}$.
 In this section we recall the most common choices and the associated structures in terms of the resulting blur operator $\tilde{X}.$

In \cite{MR1718798,MR2045058} the authors studied how the choice of the most appropriate boundary conditions can affect the precision of the image reconstruction (especially close to the boundaries) and the cost of the computation for recovering the “true” image.

In the following we consider four classical choices: the zero (Dirichlet) BCs, periodic BCs, reflective BCs and anti-reflective BCs.
\begin{enumerate}
\item Imposing zero Dirichlet BCs we assume that the signal outside the domain of the observed vector $\mathbf{g}$ is zero:
$$
\mathbf{f}_{l}=\mathbf{f}_{r}=\mathbf{0}.
$$
The linear system in (\ref{eq:syst_blur}) becomes
\begin{equation}
\label{eq:syst_blur_zero}
T \mathbf{f}=\mathbf{g},
\end{equation}
where the coefficient matrix has a Toeplitz structure. Notice that nonzero Dirichlet BCs with $\mathbf{f}_{l}$, $\mathbf{f}_{r}$ given vectors will lead to a linear system with the very same coefficient matrix as in (\ref{eq:syst_blur_zero}) but with a different righthand-side.
\item The periodic BCs require instead
$$
f_{j}=f_{n+j} \quad \text { for all } j
$$
and the linear system is
\begin{equation}
\label{eq:syst_blur_periodic}
\tilde{X}_n \mathbf{f}=\mathbf{g}, \qquad \tilde{X}_n =\left[\left(0 \mid T_{l}\right)+T+\left(T_{r} \mid 0\right)\right],
\end{equation}
where the coefficient matrix has a circulant structure since $\left(0 \mid T_{l}\right)$ and $\left(T_{r} \mid 0\right)$ are $n \times n$ Toeplitz matrices obtained by augmenting $(n-m)$ zero columns to $T_{l}$ and $T_{r}$, respectively.

\item For the Neumann reflective BCs the data outside $\mathbf{f}$ are a reflection of the data inside $\mathbf{f}$. More precisely,

$$
\left\{\begin{array} { c c c }
{ f _ { 0 } } & { = } & { f _ { 1 } } \\
{ \vdots } & { \vdots } & { \vdots } \\
{ f _ { - m + 1 } } & { = } & { f _ { m } }
\end{array} \quad \text { and } \quad \left\{\begin{array}{ccc}
f_{n+1} & = & f_{n} \\
\vdots & \vdots & \vdots \\
f_{n+m} & = & f_{n-m+1}
\end{array}\right.\right.
.$$
Thus system (\ref{eq:syst_blur}) becomes
\begin{equation}
\label{eq:syst_blur_refl}
\tilde{X_n}  \mathbf{f}=\mathbf{g},  \qquad \tilde{X_n}=\left(0 \mid T_{l}\right) Y_n+T+\left(T_{r} \mid 0\right) Y_n.
\end{equation}
The system matrix $\tilde{X_n} $ can be written as the sum of a Toeplitz matrix and a matrix
\begin{equation}\label{eq:W_R}
  W_n^{R}= \left(0 \mid T_{l}\right) Y_n+\left(T_{r} \mid 0\right) Y_n,  
\end{equation}
which is a Hankel matrix.

\item For the anti-reflective BCs, the data outside
$\mathbf{f}$ are an anti-reflection of the data inside $\mathbf{f}$.
This condition is expressed by the relations
\[
\begin{array}{cc}
{f}_{1-j}={f}_1-({f}_{j+1}-{f}_1)=2{f}_1-{f}_{j+1} & \ \ \ \mbox{for all }j=1,\ldots,m,
\\
{f}_{n+j}={f}_n-({f}_{n-j}-{f}_n)=2{f}_n-{f}_{n-j}
& \ \ \ \mbox{for all }j=1,\ldots,m.
\end{array}
\]
 If we define the vector $\mathbf{z}$, whose components
are $z_j=2\sum_{k=j}^m h_k$ for $j\le m$ and zero otherwise,
and the vector $\mathbf{w}$, whose components are
$w_{n+1-j}=2\sum_{k=j}^m h_{-k}$ for $j\le m$ and zero otherwise, then
(\ref{eq:syst_blur}) becomes
\begin{equation}\label{eq:syst_blur_antirefl}
\tilde{X_n}  \mathbf{f} = \mathbf{g}, \qquad \tilde{X_n} = \mathbf{z}e_1^T-(0 | T_l )\tilde{Y}_n + T - (T_r | 0)\hat{Y}_n+\mathbf{w}e_n^T ,
\end{equation}
where $e_k$ is the $k$th vector of the canonical basis,
$\tilde{Y}_n=\begin{bmatrix}
0 & 0\\
0 & Y_{n-1}
\end{bmatrix}$,
$\hat{Y}_n=\begin{bmatrix}
                 Y_{n-1} & 0 \\
                 0 & 0
                 \end{bmatrix}.$

In more detail, $\tilde{X_n} $ is the sum of a Toeplitz matrix and  a matrix 
\begin{equation}\label{eq:W_AR}
W_n^{AR}=\mathbf{z}e_1^T-(0 | T_l )\tilde{Y_n} - (T_r | 0)\hat{Y_n}+\mathbf{w}e_n^T,
\end{equation}
where the matrices $(0 | T_l )\tilde {Y_n}$ and
$(T_r | 0)\hat
{Y_n}$
take the following form:
\begin{equation}\label{eq:conjtilde}
\left(\begin{array}{ccccccc}
 0 & h_1 & \cdots& h_m & 0& \cdots & 0     \\
  &  \vdots &\ddots &        & &   &  \\
  &  h_m & &       & &   &      \\
 \vdots &   & & O& &   & \vdots   \\
    &   & &        & &   & \\
  & & &       & &   &    \\
 0   &  & & \cdots & &  &  0
\end{array}
\right),
\left(\begin{array}{ccccccc}
 0 &  & & \cdots & &  & 0     \\
  &   & &        & &   &  \\
  & & &       & &   &      \\
 \vdots &   & &O& &   & \vdots   \\
    &   & &        & & h_{-m}  & \\
  & & &       & \ddots & \vdots  &    \\
 0   & \cdots & 0 & h_{-m} &\cdots & h_{-1} & 0
\end{array}
\right).
\end{equation}
Then $\tilde{X}_n$ is a Toeplitz + Hankel plus 2 rank correction matrix, where the correction is placed at the first and the last columns.
\end{enumerate}

The considered BCs and related structures have a natural extension in more dimensions, for example when considering the blurring of an image in 2D.
\begin{enumerate}
    \item The zero BCs in 2 dimensions produce a block Toeplitz matrix with Toeplitz blocks. In $d$ dimensions, we have a $d$-level Toeplitz matrix.
    \item The periodic BCs in 2 dimensions produce a block circulant matrix with circulant blocks. In $d$ dimensions, we have a $d$-level circulant matrix.
    \item The extension in 2 dimensions of the Neumann reflective BCs produces a coefficient matrix $\tilde{X_n}$ which is the sum of a two-level Toeplitz matrix and $W_n^{R}$  given by the sum of
\begin{itemize}
    \item a block Toeplitz matrix with Hankel blocks;
    \item a block Hankel matrix with Toeplitz blocks;
    \item a block Hankel matrix with Hankel blocks.
\end{itemize}
In $d$ dimensions, the matrix can be described as a sum of all possible combinations of block matrices arranged on $d$ levels, where each level can independently present either a Toeplitz or Hankel structure.
    \item A natural extension of the 1D structure in (\ref{eq:syst_blur_antirefl}) is obtained imposing the following anti-reflective BCs:
$$
\begin{aligned}
f_{1-i, j} & =2 f_{1, j}-f_{i+1, j}, & & \mbox{for all } 1 \leqslant i \leqslant m, 1 \leqslant j \leqslant n, \\
f_{i, 1-j} & =2 f_{i, 1}-f_{i, j+1}, & & \mbox{for all } 1 \leqslant i \leqslant n, 1 \leqslant j \leqslant m, \\
f_{n+i, j} & =2 f_{n, j}-f_{n-i, j}, & & \mbox{for all } 1 \leqslant i \leqslant m, 1 \leqslant j \leqslant n, \\
f_{i, n+j} & =2 f_{i, n}-f_{i, n-j}, & & \mbox{for all } 1 \leqslant i \leqslant n, 1 \leqslant j \leqslant m,
\end{aligned}
$$
together with the subsequent conditions for the corners:
$$
\begin{aligned}
& f_{1-i, 1-j}=4 f_{1,1}-2 f_{1, j+1}-2 f_{i+1,1}+f_{i+1, j+1}, \\
& f_{1-i, n+j}=4 f_{1, n}-2 f_{1, n-j}-2 f_{i+1, n}+f_{i+1, n-j}, \\
& f_{n+i, 1-j}=4 f_{n, 1}-2 f_{n, j+1}-2 f_{n-i, 1}+f_{n-i, j+1}, \\
& f_{n+i, n+j}=4 f_{n, n}-2 f_{n, n-j}-2 f_{n-i, n}+f_{n-i, n-j},
\end{aligned}
$$
for $1 \leqslant i, j \leqslant m$. Note that this choice corresponds to performing anti-reflection first around the $x$ axis and then around the $y$ axis.
In the general $d$ dimensional case, the anti-reflection is done around each axis in order, such that we get a $d$ levelled matrix with a structure that is a tensorial extension of (\ref{eq:syst_blur_antirefl}). For a comprehensive discussion on the selection of this type of anti-reflection, we direct the reader to the insights provided in \cite{MR2045058}.
\end{enumerate}

\section{Main Results: the original matrix-sequences}
\label{sec:main-I}

In the case of blurring in $d$ dimensions, $d\ge 1$, ($d=1$ for signals, $d=2$ for images, $d=3$ for 3D images), given the PSF and zero-Dirichlet boundary conditions, the resulting matrix is
\[\tilde{X}_n=T_n(f), \quad f=f_{\rm  PSF}.\]
In most of the cases $f$ is a $d$ variate trigonometric polynomial or $2\pi$-periodic continuous function belonging to the Wiener class.

When considering motion blurring or shaking, the function $f$, the generating function of the Toeplitz matrix $\tilde{X}_n=T_n(f)$, is complex-valued with real Fourier coefficients. Hence we can use the Pestana-Wathen \cite{MR3323542} trick so that we solve a linear system with symmetric coefficient matrix  $X_n= Y_n\tilde{X}_n$ by using a preconditioning MINRES algorithm.

Here we refine the proposal in \cite{donatelli2022symmetrization} by analysing the case of the most precise boundary conditions, i.e. reflective \cite{MR1718798} and anti-reflective  \cite{MR2045058} BCs. As proven in \cite{Ferrari2019,MazzaPestana2021}, we remind that
\begin{equation*}
\{Y_nT_n(f)\}_n\sim_\lambda \psi_{|f|}, \ \ \ \tilde{X}_n=T_n(f)
\end{equation*}
 with $ \psi_{|f|}$ as in Definition \ref{def:psi} or any rearrangement of it (see \cite{rearr} and references therein), such as
${\rm diag}\left(|f|,-|f|)\right)$ according to (\ref{eq:psi-bis}).

\begin{Theorem}\label{th:main}
Let $A_n^{BC}$ be the matrix stemming when using a nonsymmetric PSF and either reflective or anti-reflective boundary conditions. Then,
\begin{itemize}
\item[a)]  $A_n^{BC}$ is real nonsymmetric and not diagonalizable using DCT or the DST-based anti-reflective transform \cite{AR-trasf1,AR-trasf2}.
\item[b)] Take $f$ stemming from the PSF, with $f$ being $2\pi$-periodic, continuous, complex-valued. Then $A_n^{BC}=T_n(f)+W_n^{BC}$ is such that
$\{W_n^{BC}\}_n$ is strongly clustered at zero according to Definition \ref{def-cluster} in the sense of the eigenvalues with $s=0$.
\item[c)] Take $f$ stemming from the PSF, with $f$ being $2\pi$-periodic, continuous, complex-valued. Then $A_n^{BC}=T_n(f)+W_n^{BC}$ is such that
$\{W_n^{BC}\}_n$ is strongly clustered at zero according to Definition \ref{def-cluster} in the sense of the singular values with $s=0$.
\item[d)] In the case where $f$ is only essentially bounded the matrices $A_n^{BC}=T_n(f)+W_n^{BC}$ are such that
$\{W_n^{BC}\}_n$ is weakly clustered at zero according to Definition \ref{def-cluster} in the sense of the eigenvalues with $s=0$.
\item[e)] In the case where $f$ is only essentially bounded the matrices $A_n^{BC}=T_n(f)+W_n^{BC}$ are such that
$\{W_n^{BC}\}_n$ is weakly clustered at zero according to Definition \ref{def-cluster} in the sense of the singular values with $s=0$.
\item[f)] In the case where $f$ is only essentially bounded the matrices $A_n=Y_n A_n^{BC}$ are such that $\{A_n\}_n$ satisfies the assumption of Theorem \ref{th:GLT-lambda} with $X_n=Y_n T_n(f)$, $Z_n=Y_n W_n^{BC},$ $\kappa=\psi_{|f|}$ so that
$\{A_n\}_n$ is distributed in the eigenvalue sense as $\psi_{|f|}$. Furthermore $\{A_n\}_n$ is distributed in the singular value sense as $|f|$.
\item[g)] Whenever $f$ is a $2\pi$-periodic, continuous, complex-valued function with $d=1,$ $\|Z_n\|_1=O(1)$ and hence the range of $\psi_{|f|}$ is a strong cluster for $\{A_n\}_n$ in the eigenvalue sense, while the range of ${|f|}$ is a strong cluster for $\{A_n\}_n$ in the singular value sense.
\end{itemize}
\end{Theorem}
{\bf Proof:\ }
Item a) is already known and the related analysis can be found in \cite{AR-trasf1,AR-trasf2}, where it is proven that the symmetry of the PSF is necessary and sufficient for the diagonalizable character of the resulting matrices via DCT or the DST-based anti-reflective transform, respectively.

For Item b), following \cite{fasino-tilli}, we know that the Hankel corrections in the reflective and anti-reflective settings have both trace norm bounded by a pure constant depending only on $f$, while their spectral norm is bounded by the infinity norm of $f$ (see \cite{ST-ineq}). Furthermore, in the anti-reflective case a direct check shows that also the low-rank correction have both trace norm and spectral norm bounded by absolute constants.
Hence the desired result follows from Theorem \ref{th:GLT-lambda}, with $X_n$ being the null matrices and $\kappa=0$.

For Item c), the fact that the trace norm and the low rank correction have trace norm bounded by a constant depending only on $n$ is enough to conclude that  $\{W_n^{BC}\}_n$ is strongly clustered at zero according to Definition \ref{def-cluster} in the sense of the singular values with $s=0$, taking into consideration the second item in Corollary 4.1 of \cite{taud2} with $A_n=W_n^{BC}$ and $B_n$ the null matrices.

For Item d), following \cite{fasino-tilli}, we know that the Hankel corrections in the reflective and anti-reflective settings have $o(N(\mathbf{n}))$ trace norm, while their spectral norm is bounded by a pure constant depending only on $f$ (see \cite{ST-ineq}). Furthermore, in the anti-reflective case a direct check show that also the low-rank correction has trace norm $o(N(\mathbf{n}))$ and spectral norm bounded by absolute constants.
Hence the desired result follows from Theorem \ref{th:GLT-lambda}, with $X_n$ being the null matrices and $\kappa=0$.

For Item e), the fact that the trace norm and the low rank correction have trace norm bounded by a constant depending only on $n$ is enough to conclude that  $\{W_n^{BC}\}_n$ is strongly clustered at zero according to Definition \ref{def-cluster} in the sense of the singular values with $s=0$, taking into consideration the second item in Corollary 4.1 of \cite{taud2} with $A_n=W_n^{BC}$ and $B_n$ the null matrices.

For Item f) first we employ the approximation Theorem \ref{th:GLT-lambda} and then we use the notion of a.c.s. and the approximation Theorem \ref{lem:Corollary5.1}. Indeed, by repeating the reasoning in Item d), the Hankel corrections in the reflective and anti-reflective settings have $o(N(\mathbf{n}))$ trace norm, while their spectral norm is bounded by a pure constant depending only on $f$ (see \cite{ST-ineq}). In addition, following again the proof in Item d), the low-rank correction has trace norm $o(N(\mathbf{n}))$ and spectral norm bounded by absolute constants.
Hence the desired result follows from Theorem \ref{th:GLT-lambda}, with $X_n=T_n(f)$ and $\kappa=f$.

For the singular value we use the a.c.s. notion. In fact, $\{\{Y_n T_n(f)\}_n\}_j$ is a constant a.c.s. for $ \{A_n^{BC}\}_n$. Moreover, 
$\{T_n(f)\}_n\sim_\sigma f$ and, since $Y_n$ is unitary, also $\{T_n(f)\}_n\sim_\sigma f$. Whence the use of Theorem \ref{lem:Corollary5.1} with $\phi_j=f$
so that $\phi=f$ concludes the proof. 

Finally, for Item g) it is enough to use Theorem 3.6 in \cite{MR2287378}.
\hfill $\bullet$

\begin{Remark}{\rm Few extensions of Theorem \ref{th:main})}\label{rem:extension theorem}
 Item {e)} and the second part of Item {f)} hold in the more general case of $f$ Lebesgue integrable and the related proof is contained in the subsequent Theorem \ref{th:main-ext}, using the a.c.s. machinery. 
 
 For the eigenvalues the analysis is more complicated given their sensitivity to perturbations in the non-Hermitian setting. In fact, the extension of Item d) with a Lebesgue integrable generating function, would require the extension of Theorem \ref{th:GLT-lambda} to the case of unbounded matrix-sequences, the latter being still an open problem in its maximal generality. See \cite{BaSe-NLAA} for a partial result.

\end{Remark}

\begin{Theorem}\label{th:main-ext}
With the same notation as in Theorem \ref{th:main}, with reference to  Item {e)} and the second part of Item {f)}, under the weaker assumption that the generating function $f$ is Lebesgue integrable, we have
\begin{itemize}
\item[h)] $\{W_n^{BC}\}_n\sim_\sigma 0$.
\item[i)] $\{A_n^{BC}\}_n\sim_\sigma f$. 
\item[l)] $\{Y_n A_n^{BC}\}_n\sim_\sigma f$.
\end{itemize}
\end{Theorem}
{\bf Proof:\ } 
We start proving the statement in Item h). From the structure of $W_n^{BC}$ we know that $W_n^{BC}=H_n^{BC}+R_n^{BC}$, where $H_n^{BC}$ is of $d$-level Hankel nature and $R_n^{BC}$ is zero for the reflective BCs and is of rank bounded by $C [N(\mathbf{n})]^{(d-1)/d}$ with $C$ universal constant in the case of anti-the reflective BCs. Furthermore, following \cite{fasino-tilli}, $\|H_n^{BC}\|_1=o(N(\mathbf{n}))$ in the most general case where the generating function $f$ is Lebesgue integrable. Therefore by Corollary 4.1 in \cite{taud2}, with $A_n=H_n^{BC}$, $B_n$ being the null matrix and $p=1$, we deduce that $\{W_n^{BC}\}_n$ is clustered at $s=0$ in the singular value sense.

Hence, by the definition of a.c.s. and by Theorem \ref{thm:acs_caratt}, $\{\{H_n^{BC}+R_n^{BC}\}_n\}_j$ is an a.c.s. for the null matrix-sequence whose singular distribution function is  $\phi=0$. Since $\phi_j$ does not depend on $j$ and $\phi_j$ converges to $\phi=0$ by Theorem \ref{lem:Corollary5.1} we deduce that $\{W_n^{BC}=H_n^{BC}+R_n^{BC}\}_n\sim_\sigma 0$ and Item h) is proved.

For Item i), we heavily rely on the reasoning for proving the previous statement. In fact
\[
\{A_n^{BC}\}_n=\{T_n(f)\}_n + \{W_n^{BC}\}_n.
\]
Hence $\{\{T_n(f)\}_n\}_j$ is an a.c.s. for $\{A_n^{BC}\}_n$ so that using again Theorem \ref{lem:Corollary5.1}, we deduce $\phi_j=f$ for every $j$ 
and so $\phi=f$ for which 
\[
\{A_n^{BC}\}_n\sim_\sigma f. 
\]
Finally Item l) is a direct consequence of Item i), taking into account that $Y_n$ is unitary and hence, independently of the matrix-size,
the matrices $A_n^{BC}$ and $Y_n A_n^{BC}$ share the same singular values.

\hfill $\bullet$

\begin{Remark}{\rm (The GLT machinery)}\label{rem:glt}
As a counterpart of Theorem \ref{th:main-ext}, the GLT analysis would have lead to $\{W_n^{BC}\}_n\sim_{\rm GLT} 0$, $\{A_n^{BC}\}_n\sim_{\rm GLT}  f$
$\{Y_n A_n^{BC}\}_n\sim_\sigma f$. This machinery is indeed necessary when proving results for spatially varying blurring analogous to those previously shown, where the Toeplitz structure is algebraically lost and survives asymptotically only at a local level.
\end{Remark}

\begin{Remark}{\rm (The preconditioned matrix-sequences)}\label{rem:prec}
Taking inspiration from Theorem \ref{th:main-ext}, as already observed in Remark \ref{rem:glt}, the GLT analysis could be employed for alternative proofs. However, the main point is that varying $m$ we have several GLT $*$-algebras and the power of the algebraic structure is that any algebraic operation on GLT matrix-sequences reflects analogously on the same algebraic operation on the associated GLT symbols. The latter has a direct  use in the preconditioning spectral/singular value analysis.
Indeed, for any matrix-sequence $\{P_n \}_n\sim_{\rm GLT} f$ with $P_n$ invertible we have $\{P_n^{-1} A_n^{BC}\}_n\sim_{\rm GLT} 1$ so that
$\{P_n^{-1} A_n^{BC}\}_n\sim_\sigma 1$ and a weak singular value clustering occurs: here $P_n$ can be chosen in the circulant, $\tau$ algebra or even in the algebra diagonalized by the anti-reflective algebra. However, according to the GLT theory, the latter is true only when the generating function $f$ is sparsely vanishing i.e. when the sets of its zeros has zero Lebesgue measure (see \cite{GSI}[Chapter 9, Item GLT5, p. 170] for $d=1$ and \cite{GSII}[Chapter 6, Item GLT5, p. 118] for $d\ge 2$). 

Regarding the spectral analysis, the situation is more involved but we could rely on the advanced tools in \cite{BaSe-NLAA}.

However, in the present work we limit ourselves to the nonpreconditioned case and we leave the analysis in the preconditioning setting for a future research step, taking into account regularizing preconditioners and the fact that in practice it may happen that the GLT symbol $f$ fails to be sparsely vanishing as it is clear for instance in Figure~\ref{fig:comparison_eig_R} and in Figure~\ref{fig:comparison_eig_AR}. 
\end{Remark}

\section{Numerical Experiments}\label{sec:num}
In the current section we perform several experiments to numerically confirm the distribution results of Theorem \ref{th:main} and their computational implications. In particular, we focus on the two-dimensional case, dealing with image reconstruction by means of the GMRES method.

While it has been established that the MINRES algorithm is a regularization method for symmetric ill-posed linear systems \cite{MR1413298}, the same cannot be generally said for the GMRES algorithm when dealing with non-Hermitian matrices. When the matrix deviates significantly from a Hermitian structure, the GMRES method may struggle in accurately approximating the solution to an ill-posed problem \cite{MR1885302}. The presence of non-real eigenvalues in a matrix is a clear indication of its non-Hermitian nature. The more non-real eigenvalues a matrix has and the larger the imaginary parts of such eigenvalues are, the more it deviates from the real Hermitian subspace. 

Building upon the latter informal statement and to further substantiate this line of reasoning, in Subsection \ref{ssec:eigenvalue_plots} we study the complex eigenvalues of the flipped system matrix in the case of a non-symmetric blur, imposing both reflective and anti-reflective BCs. Moreover, for the same test problem, we numerically investigate Item f) of Theorem \ref{th:main} comparing the eigenvalues of the flipped blurring matrix to a uniform sampling of the function $\psi_{|f|}$.

In Subsection \ref{ssec:image_reconstruction} we consider two different examples of restoration of images contaminated by blur and noise.  In both cases we show that the flipping matrix $Y_n$ is indeed an efficient regularising preconditioner for the application of the GMRES method to the considered linear systems. Furthermore, we comment on the comparison between the reconstruction quality resulting from all the boundary conditions presented in the previous sections.

\subsection{Eigenvalue Plots}\label{ssec:eigenvalue_plots}
In the present subsection, we study the eigenvalues associated with flipped system matrices, focusing on a scenario involving a non-symmetric blur. In particular, in the whole subsection we consider the speckle blur depicted in Figure~\ref{fig:squirrel_images} and properly cropped to match with the size of the blurring matrices that we consider for an image of size $64 \times 64$.

In scenarios involving periodic and zero boundary conditions, the matrices are circulant and Toeplitz, respectively. Previous studies, such as \cite{Ferrari2021}, have numerically illustrated the eigenvalue distributions for these flipped matrix-sequences. We remark again that in these cases the flipped matrices are Hermitian, resulting in real eigenvalues. The latter feature does not generally apply to reflective and anti-reflective BCs. Therefore, our focus will be on these latter two cases.

\begin{figure}[htbp]
     \includegraphics[width=\textwidth]{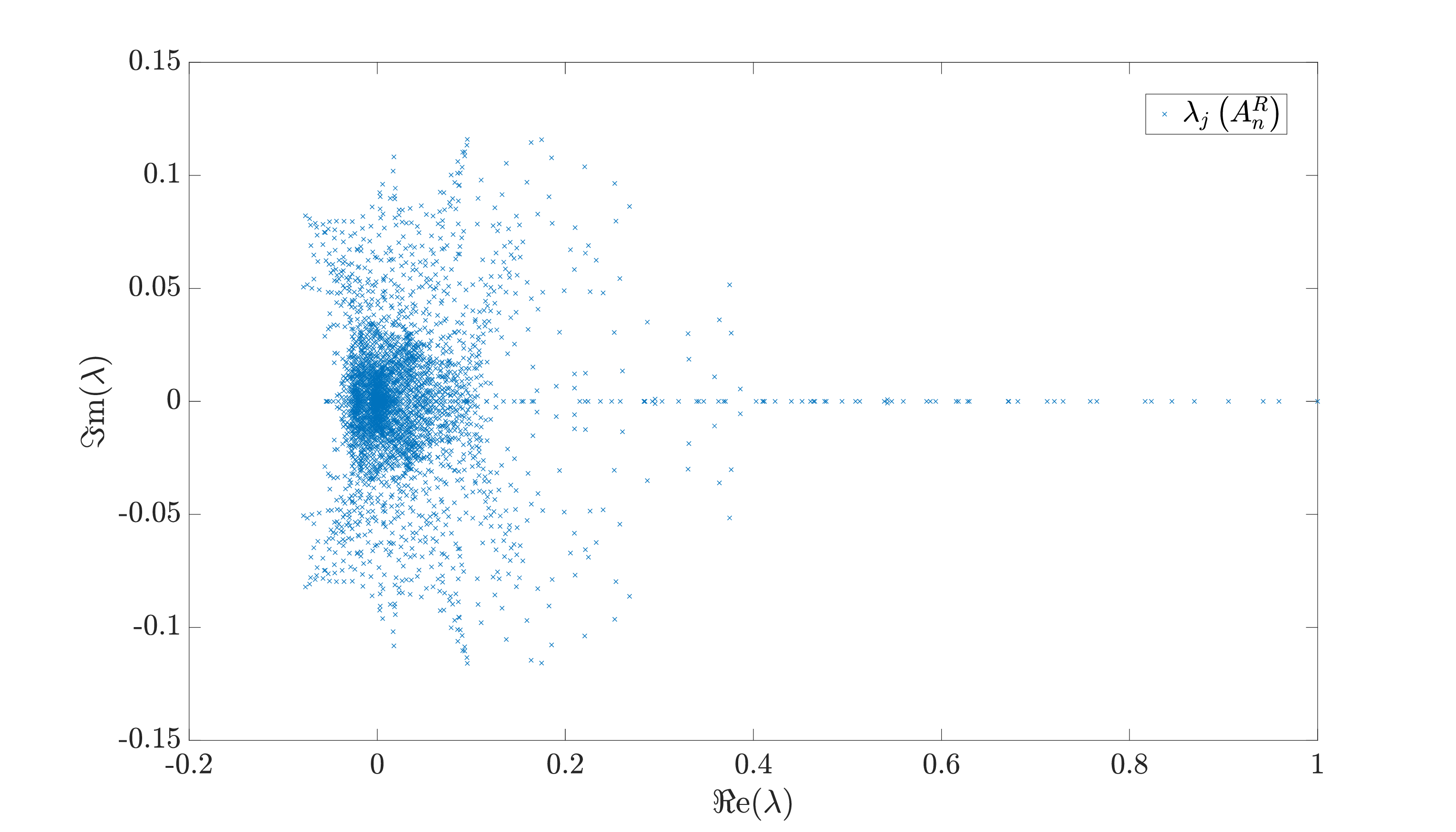}
    \caption{Eigenvalues in the complex plane of the non-flipped matrix $A_n^{R}$ of size $4096 \times 4096$.   }
    \label{fig:complex_eig_R}
\end{figure}

\begin{figure}[htbp]
     \includegraphics[width=\textwidth]{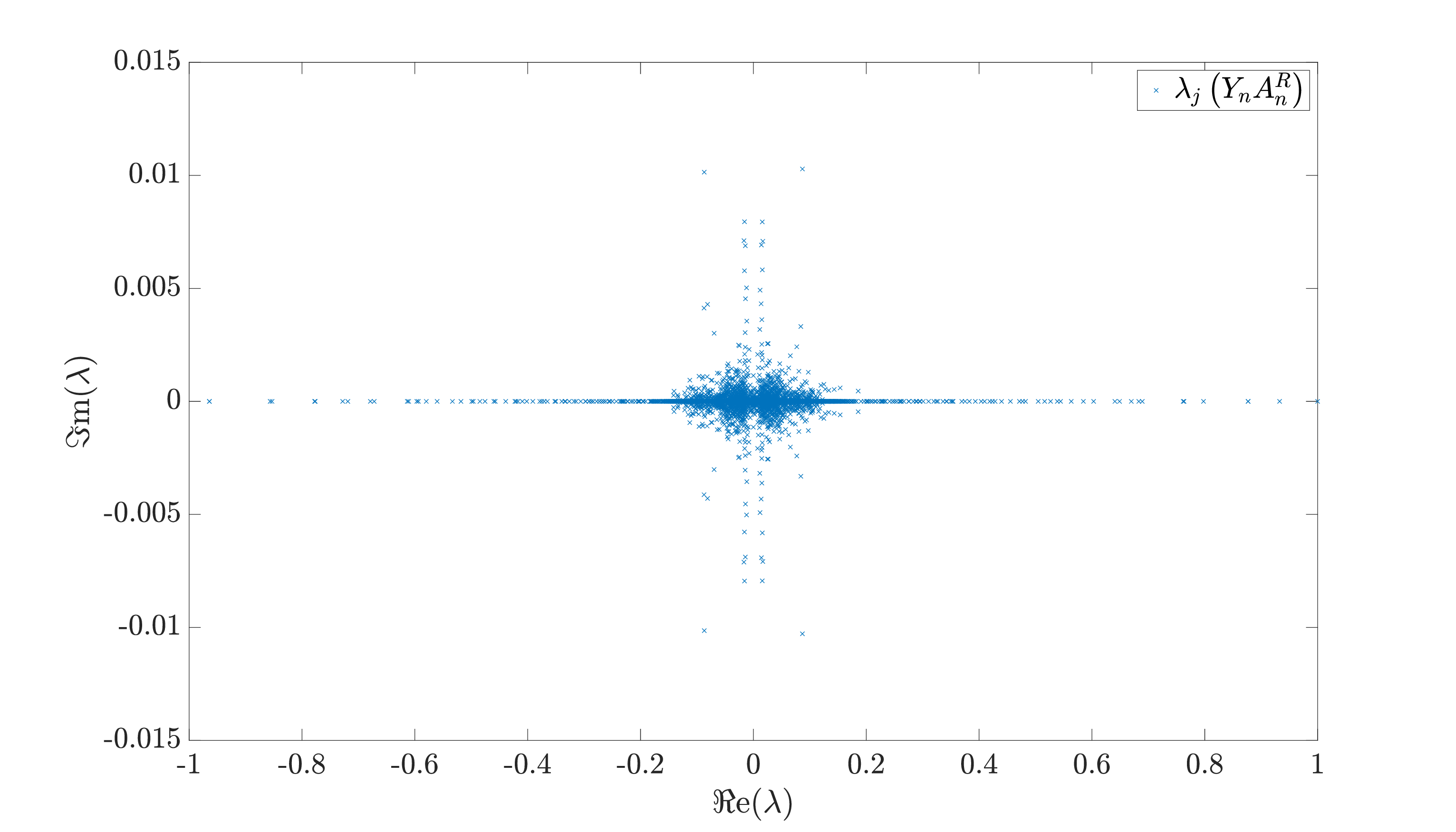}
    \caption{Eigenvalues in the complex plane of the flipped matrix $Y_n A_n^{R}$ of size $4096 \times 4096$.   }
    \label{fig:complex_eig_YR}
\end{figure}

In Figures \ref{fig:complex_eig_R}--\ref{fig:complex_eig_YR} we consider reflective BCs and we plot the eigenvalues of the non-flipped matrix $A_n^{R}$ (Figure \ref{fig:complex_eig_R}) and of the flipped matrix $Y_n A_n^{R}$ (Figure \ref{fig:complex_eig_YR}) in the complex plane. We see that the flip operation notably influences the eigenvalues, resulting in a higher proportion of real values and reducing by an order of magnitude their imaginary parts. For example, the $4096 \times 4096$ non-flipped matrix has 1939 eigenvalues with a numerically non-zero imaginary part and the flipped matrix has 1289 eigenvalues with a numerically non-zero imaginary part. 

\begin{figure}[htbp]
     \includegraphics[width=\textwidth]{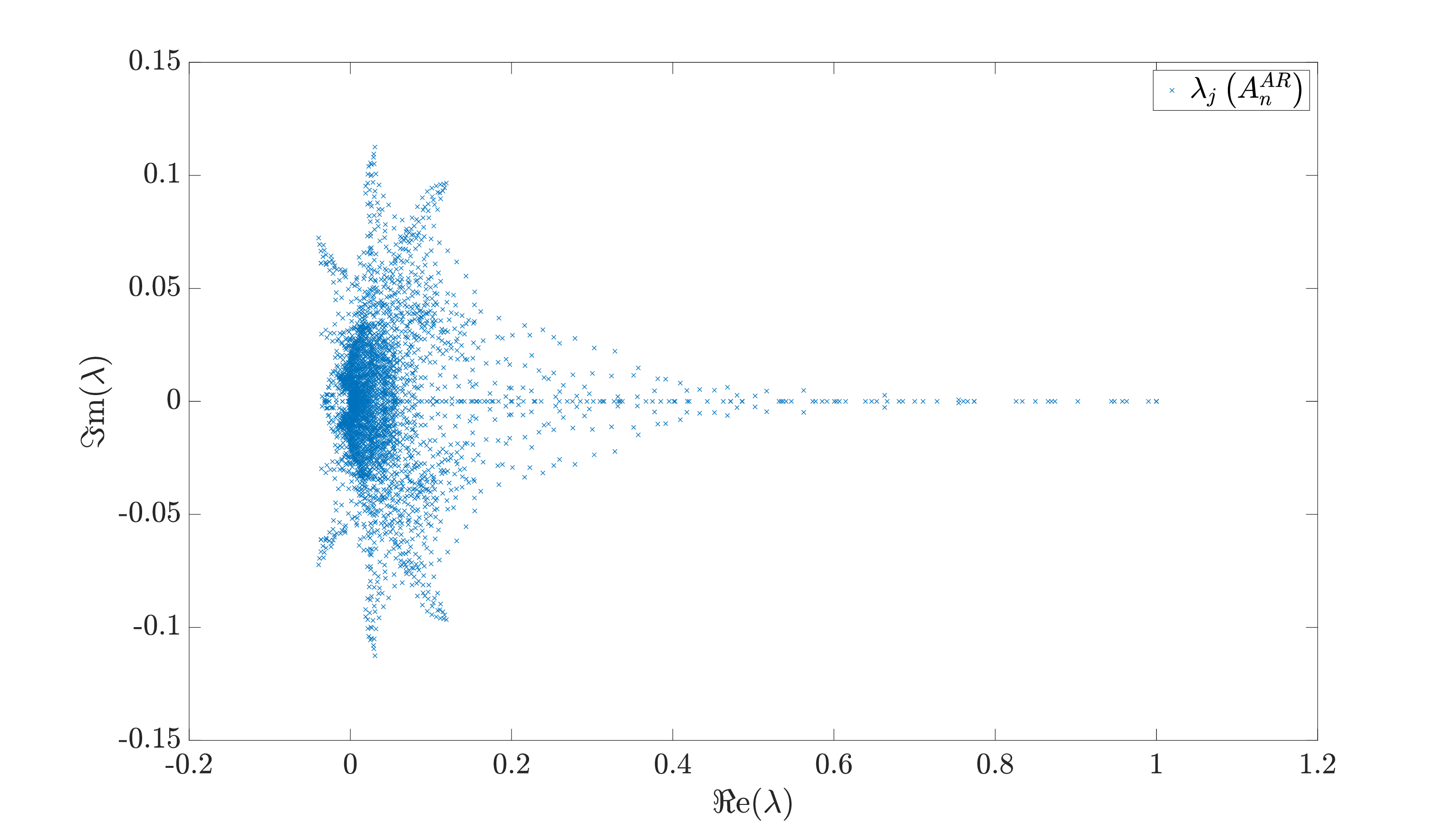}
    \caption{Eigenvalues in the complex plane of the non-flipped matrix $A_n^{AR}$ of size $4096 \times 4096$.   }
    \label{fig:complex_eig_AR}
\end{figure}

\begin{figure}[htbp]
     \includegraphics[width=\textwidth]{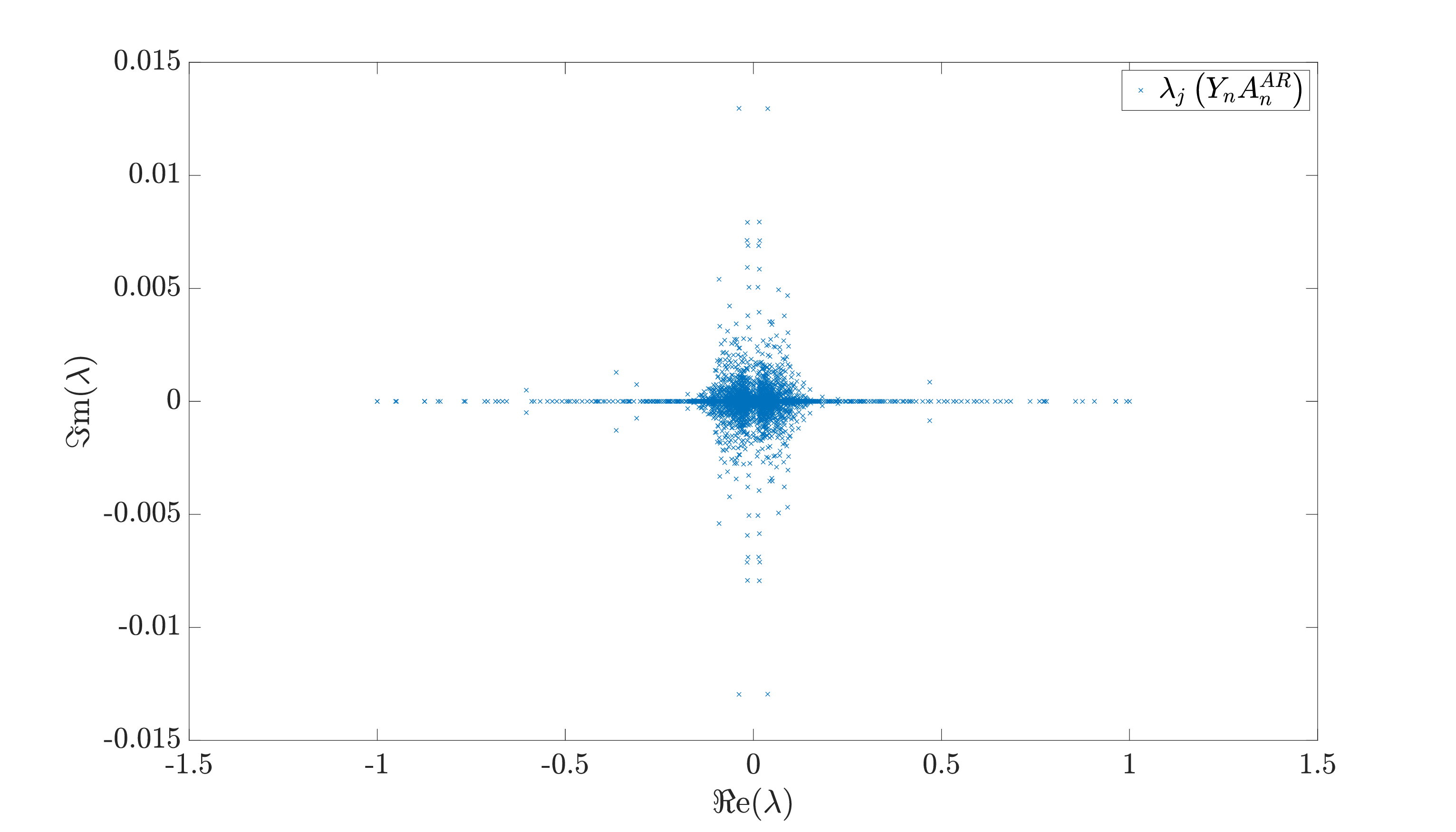}
    \caption{Eigenvalues in the complex plane of the flipped matrix $A_n^{AR}$ of size $4096 \times 4096$.   }
    \label{fig:complex_eig_YAR}
\end{figure}

Analogously, we analyze the case of anti-reflective BCs plotting the eigenvalues of $A_n^{AR}$ in Figure \ref{fig:complex_eig_AR} and of $Y_n A_n^{AR}$ in Figure \ref{fig:complex_eig_YAR}. Also in this case, the number of real eigenvalues and the magnitude of their imaginary parts decreases after performing the flip operation. The non-flipped matrix has 1941 eigenvalues with a numerically non-zero imaginary part and the flipped matrix has 1378 eigenvalues with a numerically non-zero imaginary part. We highlight that in the anti-reflective scenario the reduction in the number of non-real eigenvalues is proportionally less evident with respect to the reflective case.

\begin{figure}[htbp]
     \includegraphics[width=\textwidth]{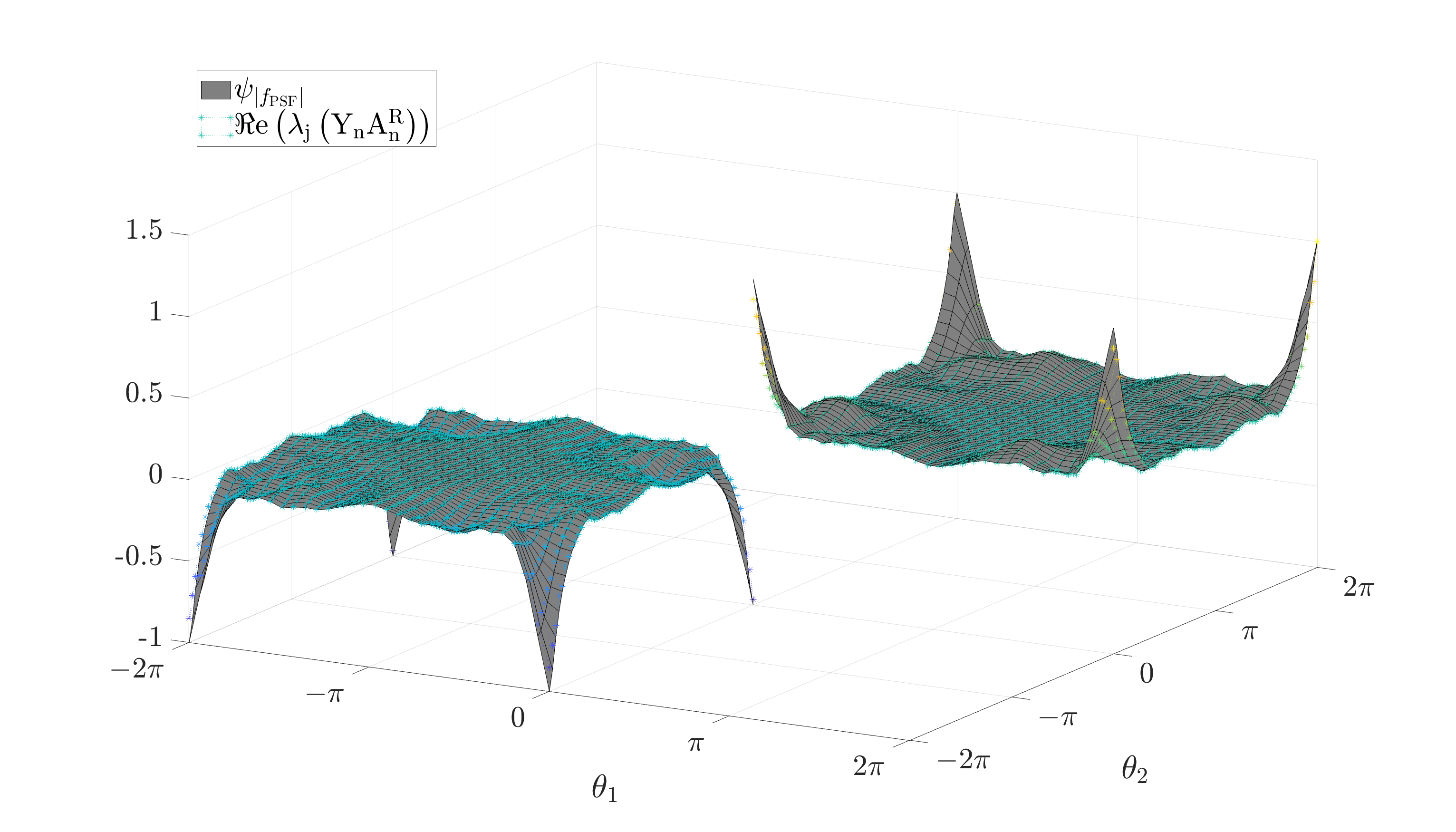}
    \caption{Comparison of the real parts of the eigenvalues of the flipped matrix $Y_nA_n^{R} $ of size $4096 \times 4096$ with a uniform sampling of the function $\psi_{|f_{PSF}|}$. }
    \label{fig:comparison_eig_R}
\end{figure}

\begin{figure}[htbp]
     \includegraphics[width=\textwidth]{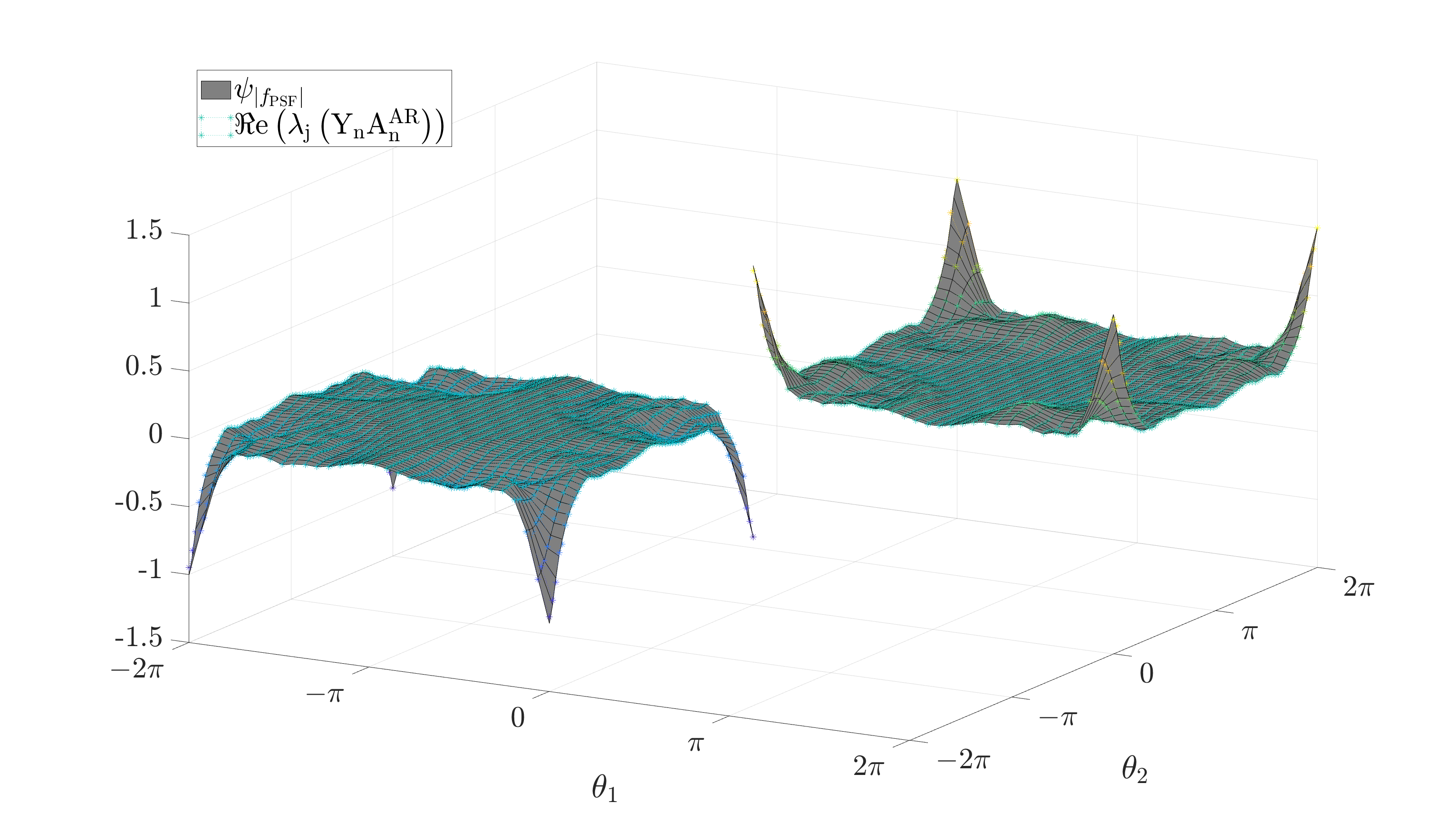}
    \caption{Comparison of the real parts of the eigenvalues of the flipped matrix $Y_nA_n^{AR} $ of size $4096 \times 4096$ with a uniform sampling of the function $\psi_{|f_{PSF}|}$.}
    \label{fig:comparison_eig_AR}
\end{figure}

The final objective of this section is to carry out a numerical validation of Item f) from Theorem \ref{th:main}, where we compare the eigenvalues of the flipped blurring matrices with a uniform sampling of the function $\psi_{|f|}$. Given that the eigenvalues of these matrices are complex, a direct visual comparison is difficult. Therefore, we focus on the real parts of the eigenvalues. For insights into the behavior of the imaginary parts, we refer readers to Figures \ref{fig:complex_eig_R}--\ref{fig:complex_eig_YAR}, where we can observe that the largest imaginary parts move from around $10^{-1}$ to around $10^{-2}$. In this context, the complex eigenvalues tend to have negligible imaginary parts as the matrix-size increases. Furthermore, the real parts are approximated by the spectral symbol, as expected, which is real by definition and in addition no outliers seem to be present. The comparison results are depicted in Figure \ref{fig:comparison_eig_R} for reflective boundary conditions and in Figure \ref{fig:comparison_eig_AR} for anti-reflective boundary conditions. In both instances we observe that the approximation of the eigenvalues by a sampling of $\psi_{|f|}$ implied by the asymptotic distribution result can already be appreciated for a relatively small matrix size.

\subsection{Image Reconstruction}\label{ssec:image_reconstruction}

The current section focuses on the problem of restoration of images contaminated by blur and noise. The blurring operator $A$ is represented by a severely ill-conditioned matrix, as we explained in the theoretical sections {and as we numerically confirmed in Figures \ref{fig:comparison_eig_R}--\ref{fig:comparison_eig_AR}, where it can be observed that a whole surface of eigenvalues is close to 0}. The reconstruction problem consists in finding a numerical solution of the linear system $A\mathbf{f} = \tilde{\bg}$ as accurately as possible. In real-world applications, the observed image is also contaminated by noise. To take this into consideration, in all examples we change the right-hand side as follows:
 \[\tilde{\bg}=A\mathbf{f} +\frac{\bzeta}{\|\bzeta\|}\gamma\|A\mathbf{f}\|,\] 
	where $\mathbf{f}$ is the exact image, $\bzeta$ 
 is a random Gaussian white noise vector and $\gamma$ specifies the noise level. 
 
In general, we do not expect convergence of a method like the GMRES applied to the linear system with the altered right-hand side $\tilde{\bg}$, because $\tilde{\bg}$ might not be in the numerical image of the ill-conditioned matrix $A$. Nevertheless, our goal is to observe a semi-convergence behavior in the GMRES method, indicating partial progress towards a solution. In order to quantify the performance of the GMRES as a restoration method, we first consider the Relative Restoration Error $ \text{RRE}(\tilde{\mathbf{f}}) = \frac{\|\tilde{{\mathbf{f}}} - {\mathbf{f}}\|_2}{\|{\mathbf{f}}\|_2}$, where $\tilde{{\mathbf{f}}}$ is the reconstructed image and ${\mathbf{f}}$ is the original image. Looking at the error curves, the best reconstruction is the one with the lowest RRE value.

We also consider a metric that gives information on the visual image reconstruction quality, that is, Peak Signal-to-Noise Ratio (PSNR). PSNR compares the maximum pixel value of the image $\max(\mathbf{f})$ to the noise as
\[ \text{PSNR}(\tilde{\mathbf{f}}) = 20 \cdot \log_{10}\left(\frac{n^2\max(\mathbf{f})}{\|\tilde{\mathbf{f}}-\mathbf{f}\|_2}\right), \]
that is, the higher the PSNR value, the better the deblurring performance.

In what follows, we describe two separate restoration examples, where the implementation of both the blurring operators and the GMRES method make use of \IRTools \cite{ir}. In evaluating an appropriate stopping criterion, we test the discrepancy principle, which halts the iterative process once the residual vector norm falls below the noise level.

\begin{figure}[htbp]
		 \includegraphics[width=0.9\textwidth]{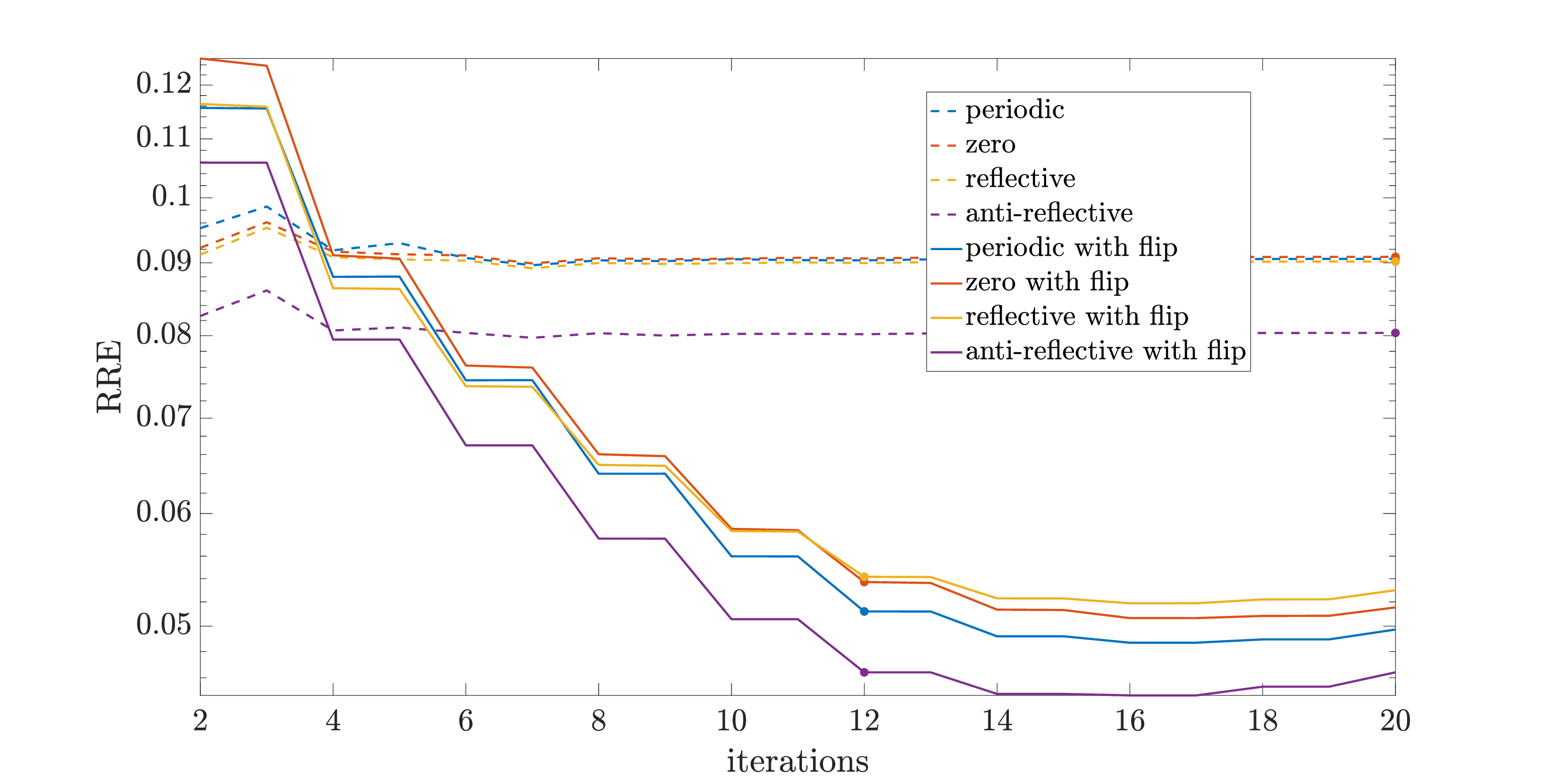}
		\caption{Flowers image reconstruction: error and iteration behaviours  of the GMRES method applied to the non-flipped (dashed lines) and flipped (continuous lines) system with different boundary conditions. The dots on the curves mark the iterations where the discrepancy principle is met.}
		\label{fig:flowers_errors}
	\end{figure}

\begin{figure}[htbp]
		 \includegraphics[width=0.3\textwidth]{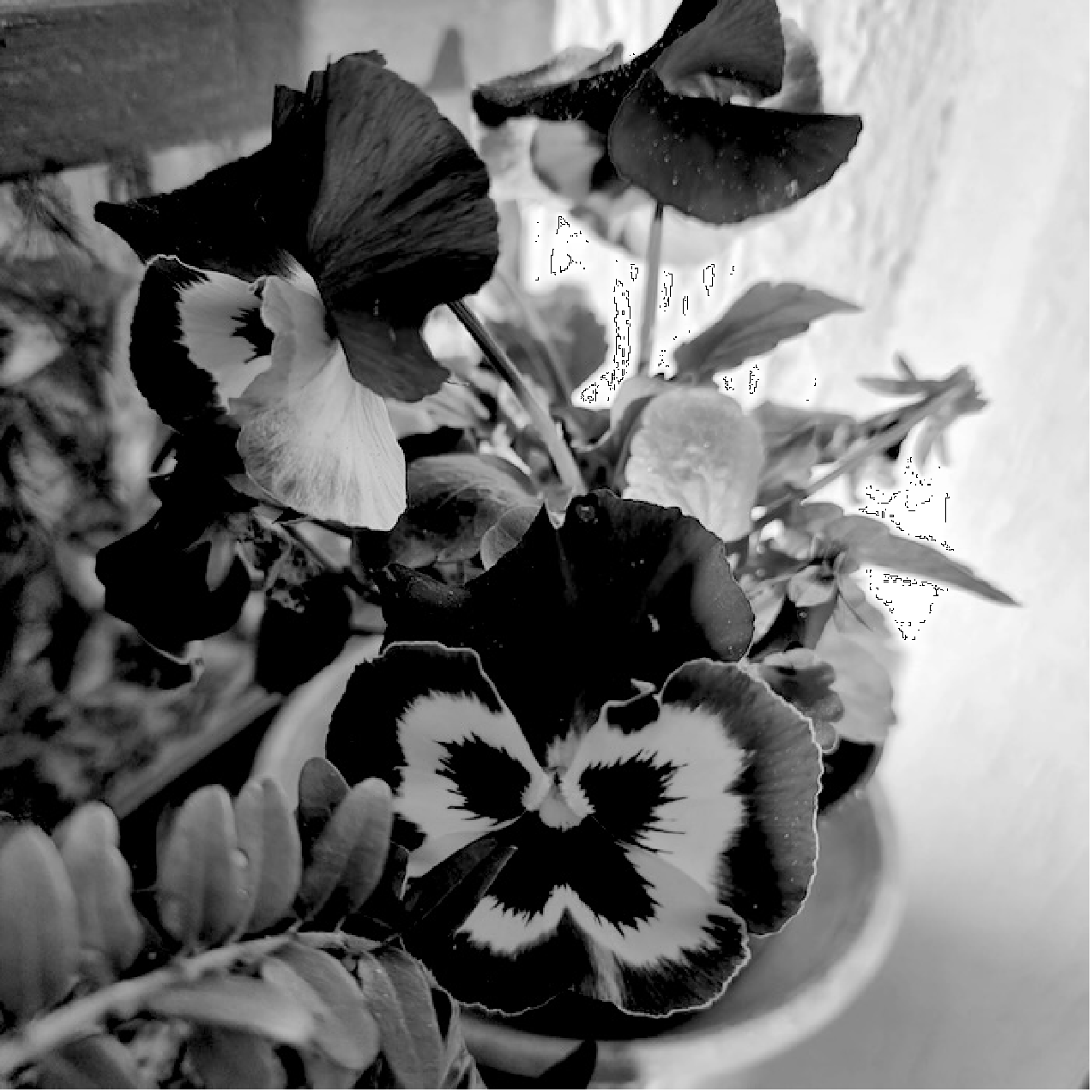}
		 \includegraphics[width=0.3\textwidth]{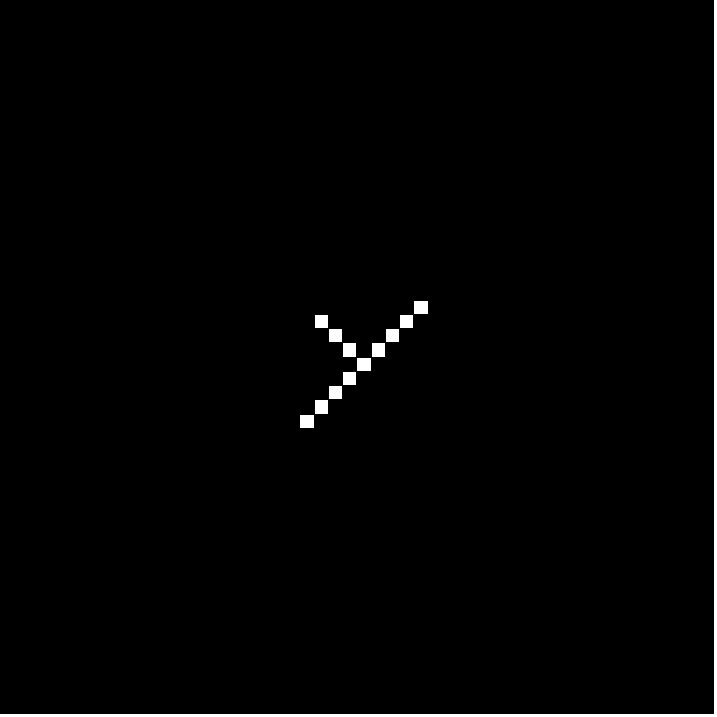}
		 \includegraphics[width=0.3\textwidth]{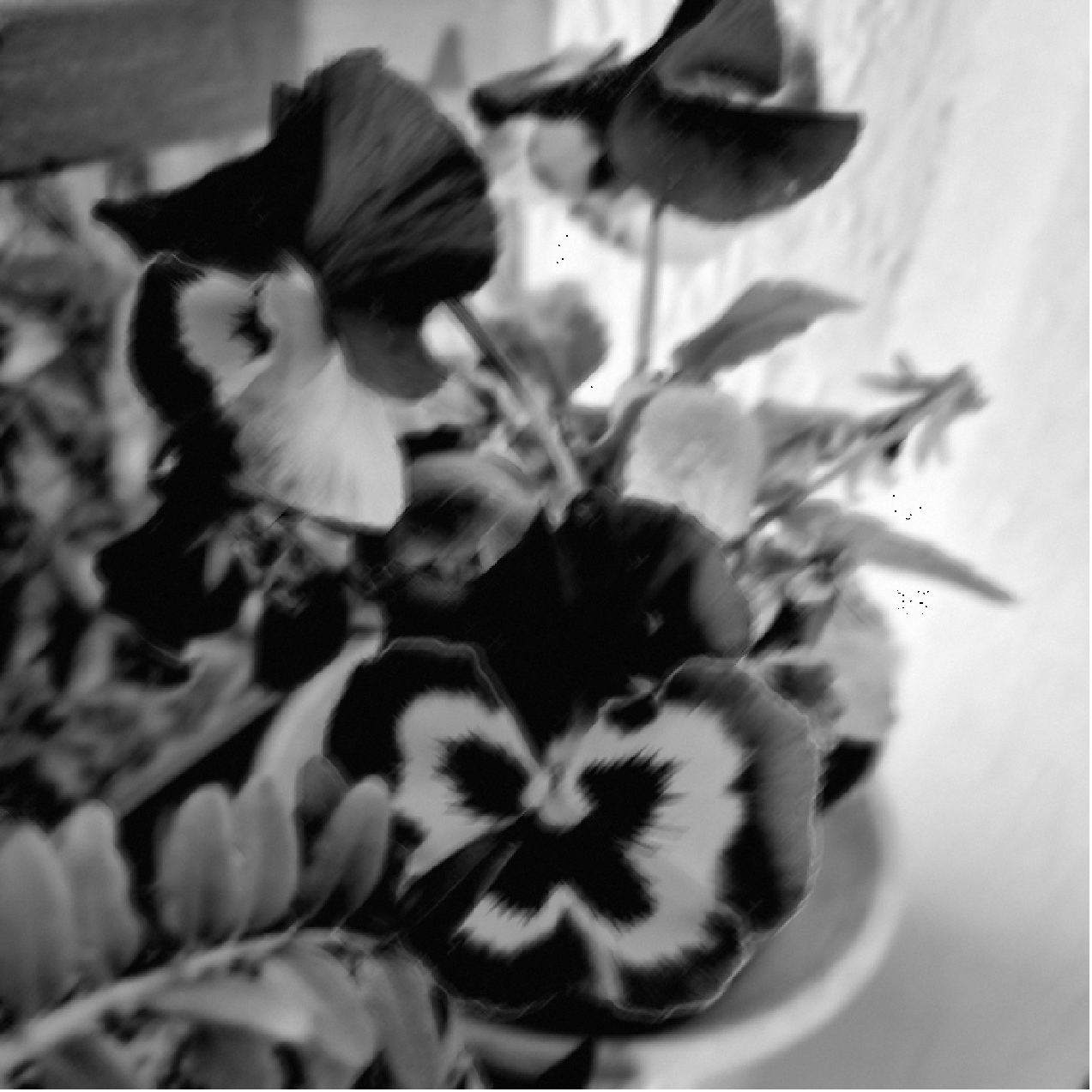}
		\caption{Flowers image, from left to right: the exact image, the motion PSF cropped to $[h_{i,j}]_{i,j=-25}^{25}$, and blurred and noisy image.}
		\label{fig:flowers_images}
	\end{figure}

\paragraph{Flowers Example}
	The first problem that we consider is the reconstruction of an image of size $700\times 700$ pixels representing flowers. The image is contaminated with a motion blur in two directions and a noise of level $\gamma=0.01$. The exact image, the PSF, and the blurred and noisy image are shown in Figure~\ref{fig:flowers_images}.
 In Figure \ref{fig:flowers_errors}, the error patterns of the GMRES method for image reconstruction under various boundary conditions (outlined in Section \ref{sec:BCs}) are presented. The dots on the curves mark the iterations where the discrepancy principle is met. A dot on the final iteration implies that the discrepancy principle does not stop the iterative process.
A first observation is that the flipping strategy is a regularising preconditioner also in the case of reflective and anti-reflective BCs, as Theorem \ref{th:main} suggests. Indeed, in all cases the continuous lines are significantly lower than the dashed ones, meaning that flipping improves the quality of the reconstruction.  Moreover, in the non-flipped case the discrepancy principle fails, while in the flipped scenario it stops the iteration when the error is close to the lowest possible.
We also remark that in this example the anti-reflecting BCs are the most suitable ones since the purple dashed line is lower than the others. The improvement is maintained also after the flip operation.

\begin{figure}[htbp]
		 \includegraphics[width=0.9\textwidth]{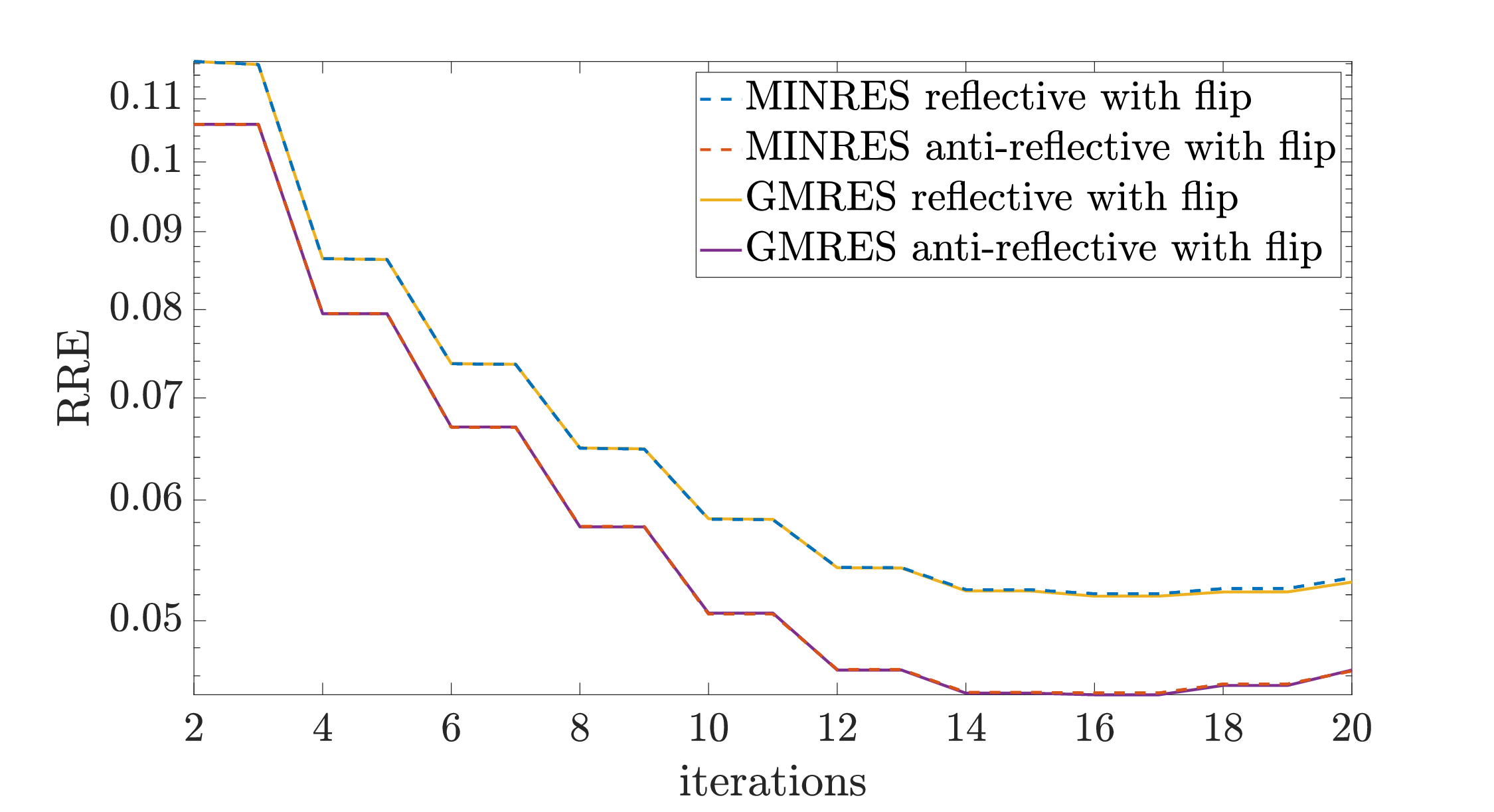}
		\caption{Flowers image reconstruction: error and iteration behaviours  of the GMRES method (continuous lines) and MINRES method (dashed lines) with different reflective and anti-reflective BCs.}
		\label{fig:flowers_errors_minres}
	\end{figure}

\begin{figure}[htbp]
		 \includegraphics[width=0.3\textwidth]{flowers_image_true.eps}
		 \includegraphics[width=0.3\textwidth]{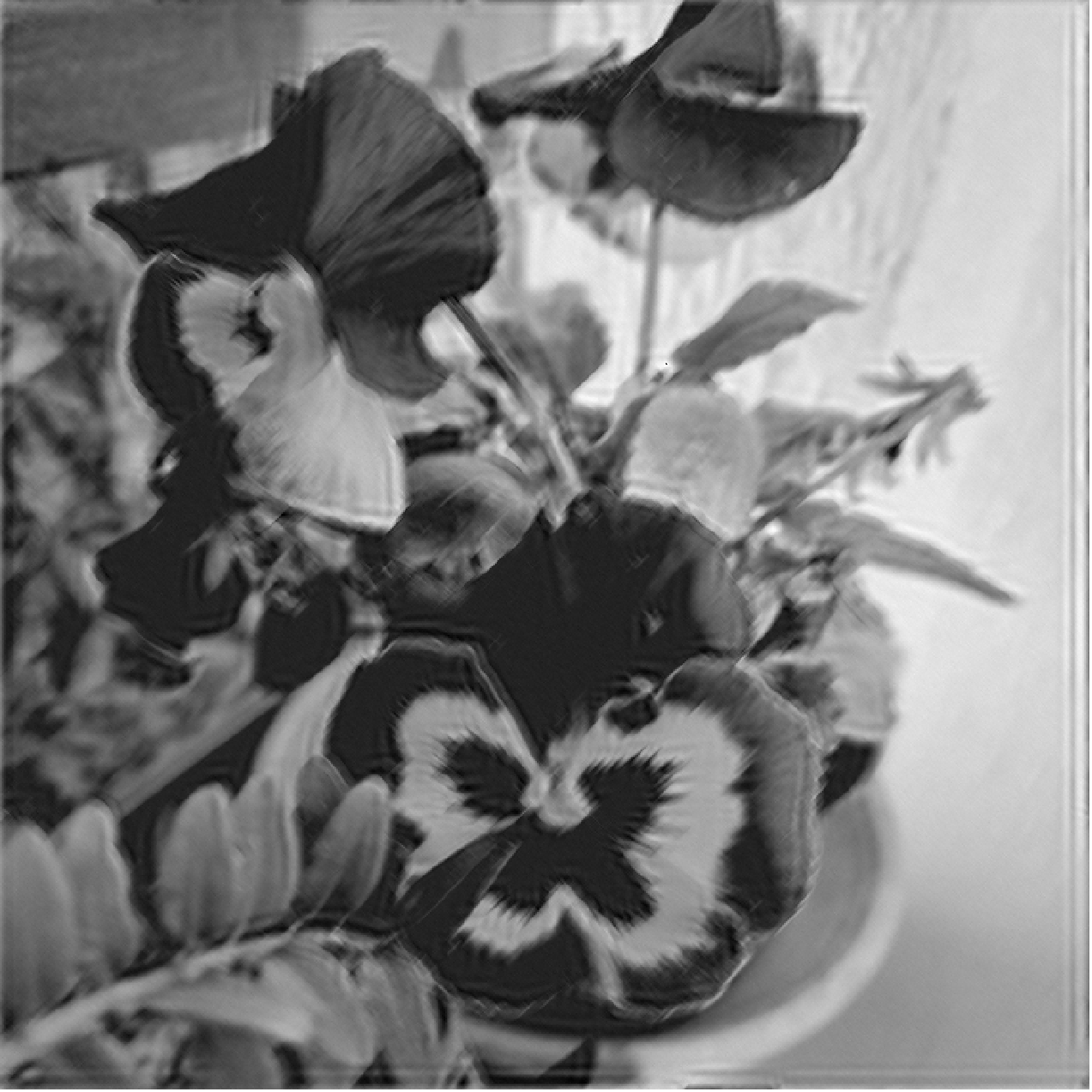}
		 \includegraphics[width=0.3\textwidth]{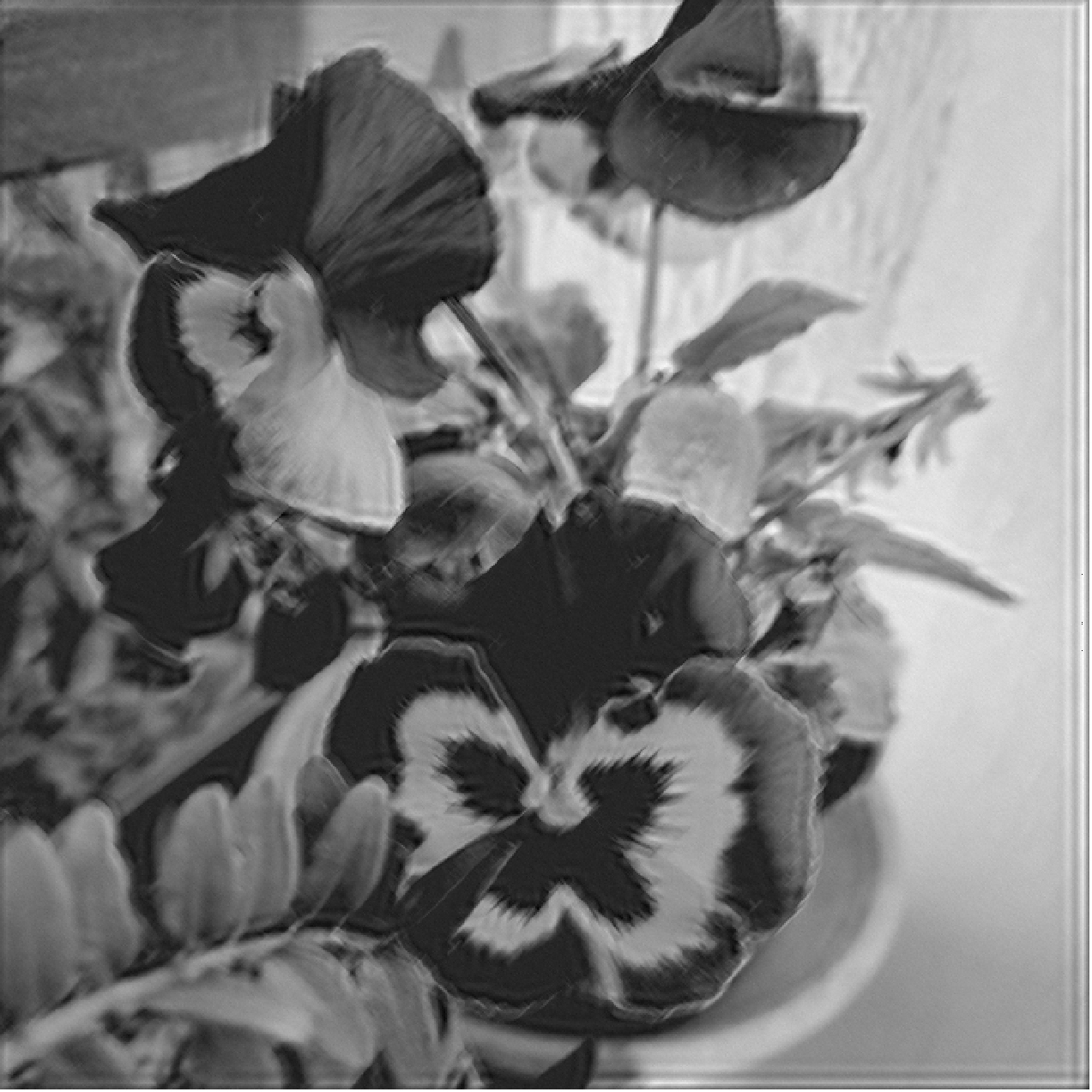}\\
   \includegraphics[width=0.3\textwidth]{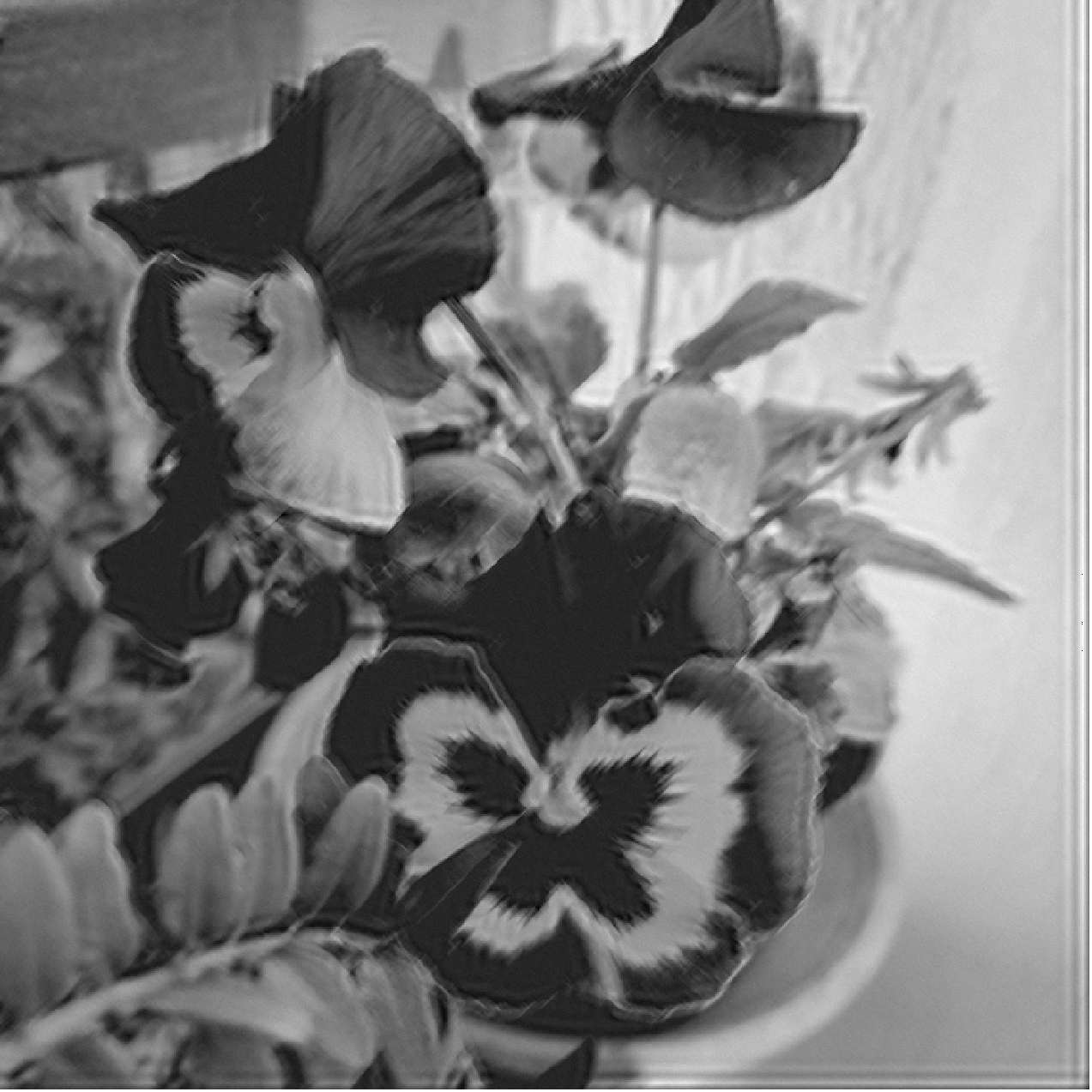}
   \includegraphics[width=0.3\textwidth]{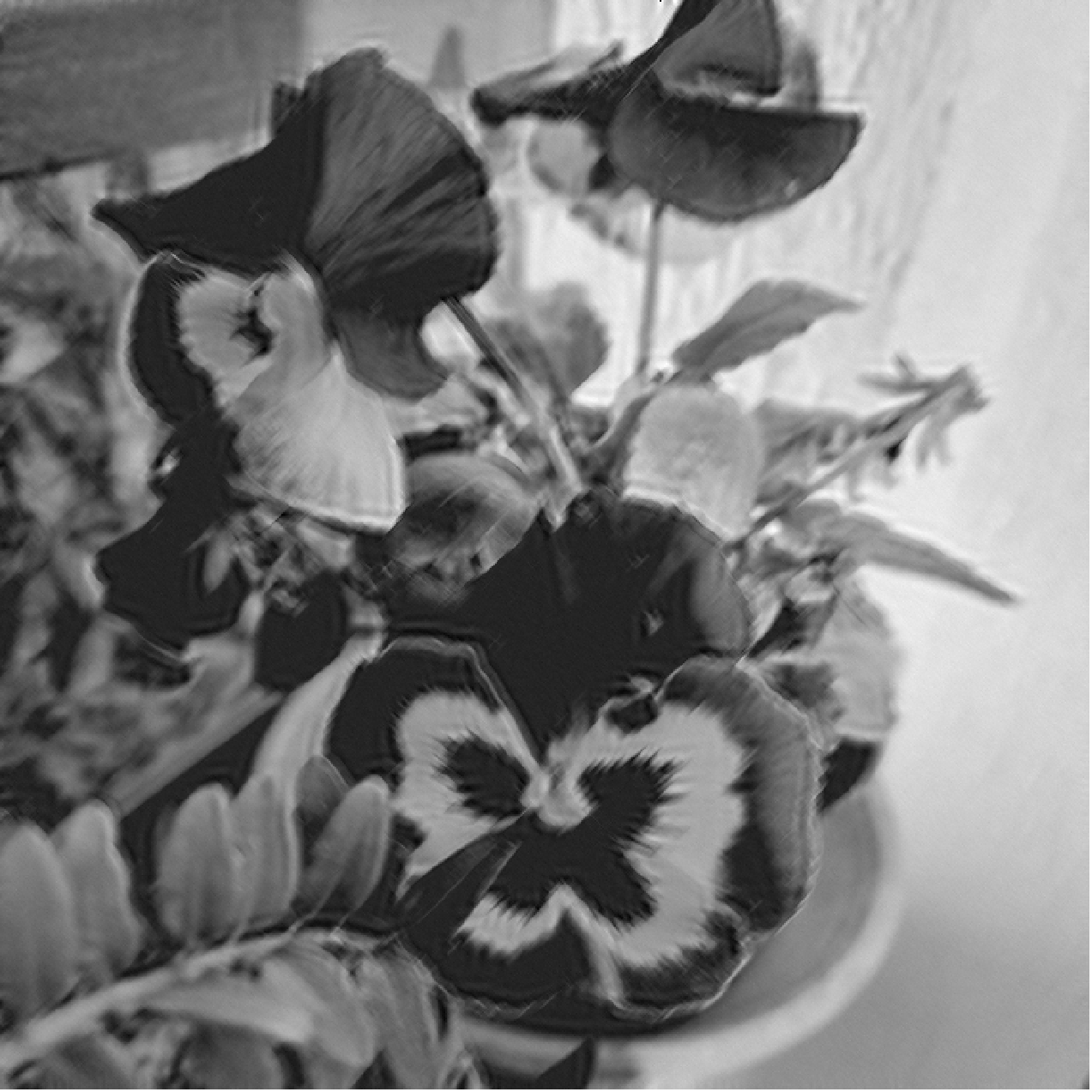}
  		 \includegraphics[width=0.3\textwidth]{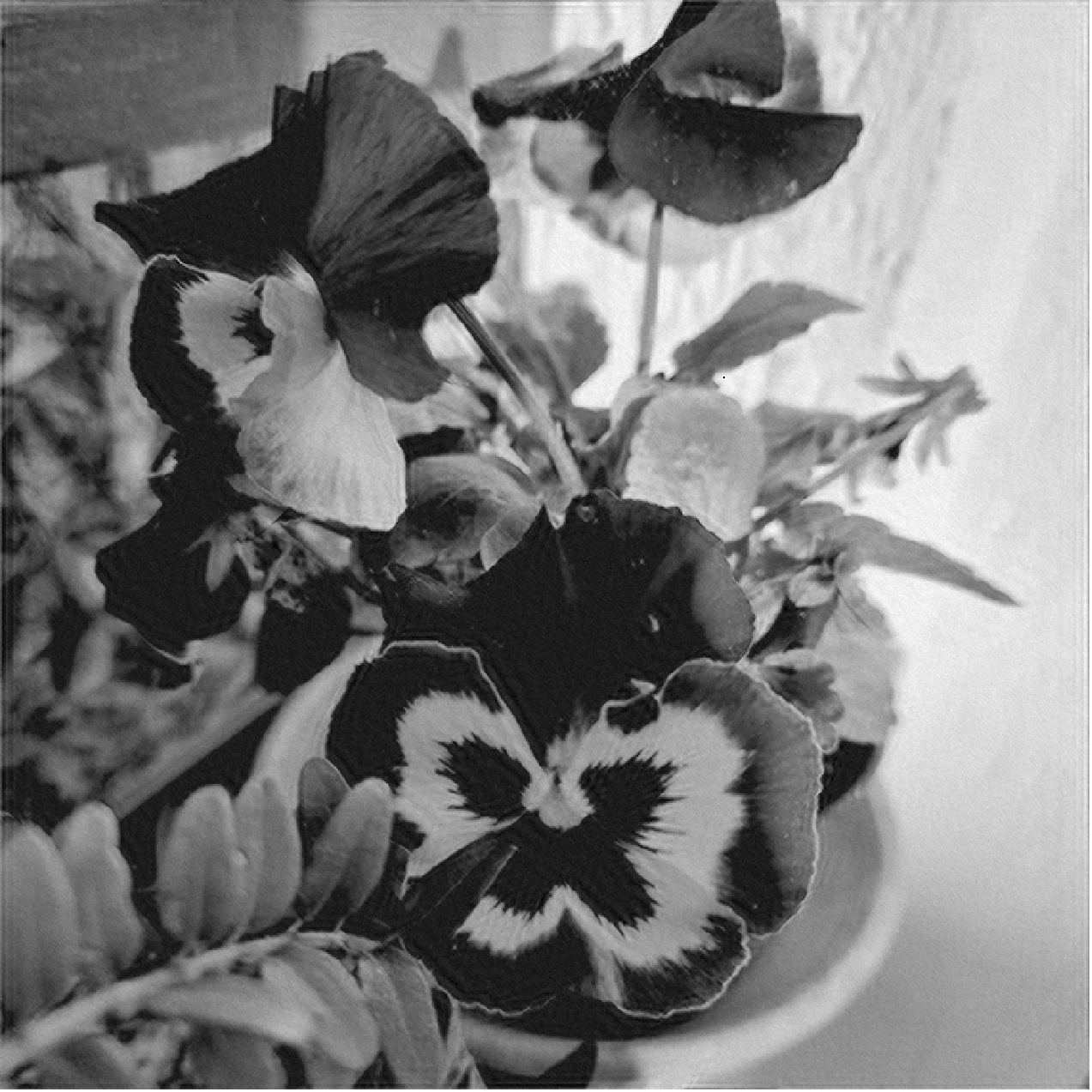}\\
		 \includegraphics[width=0.3\textwidth]{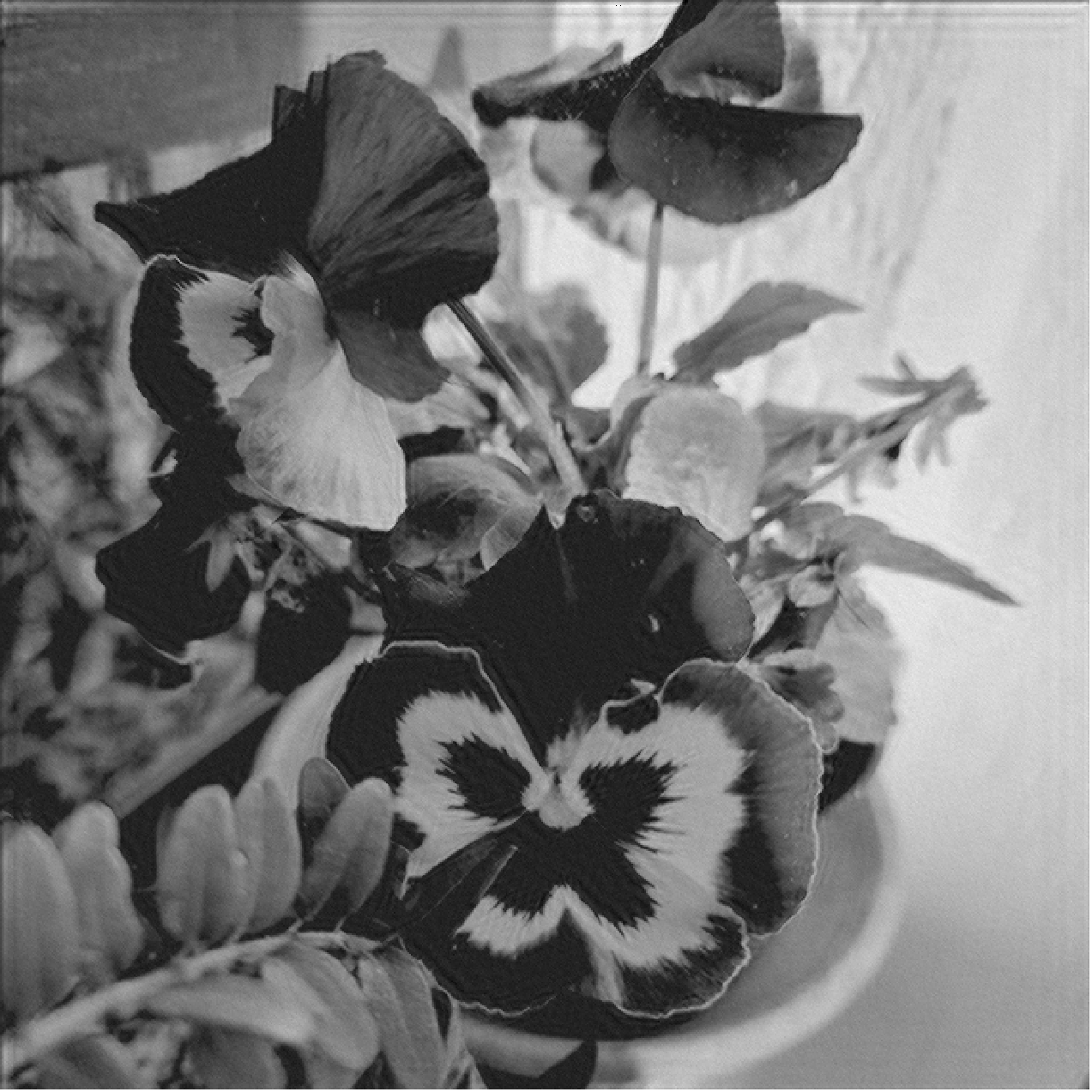}
   \includegraphics[width=0.3\textwidth]{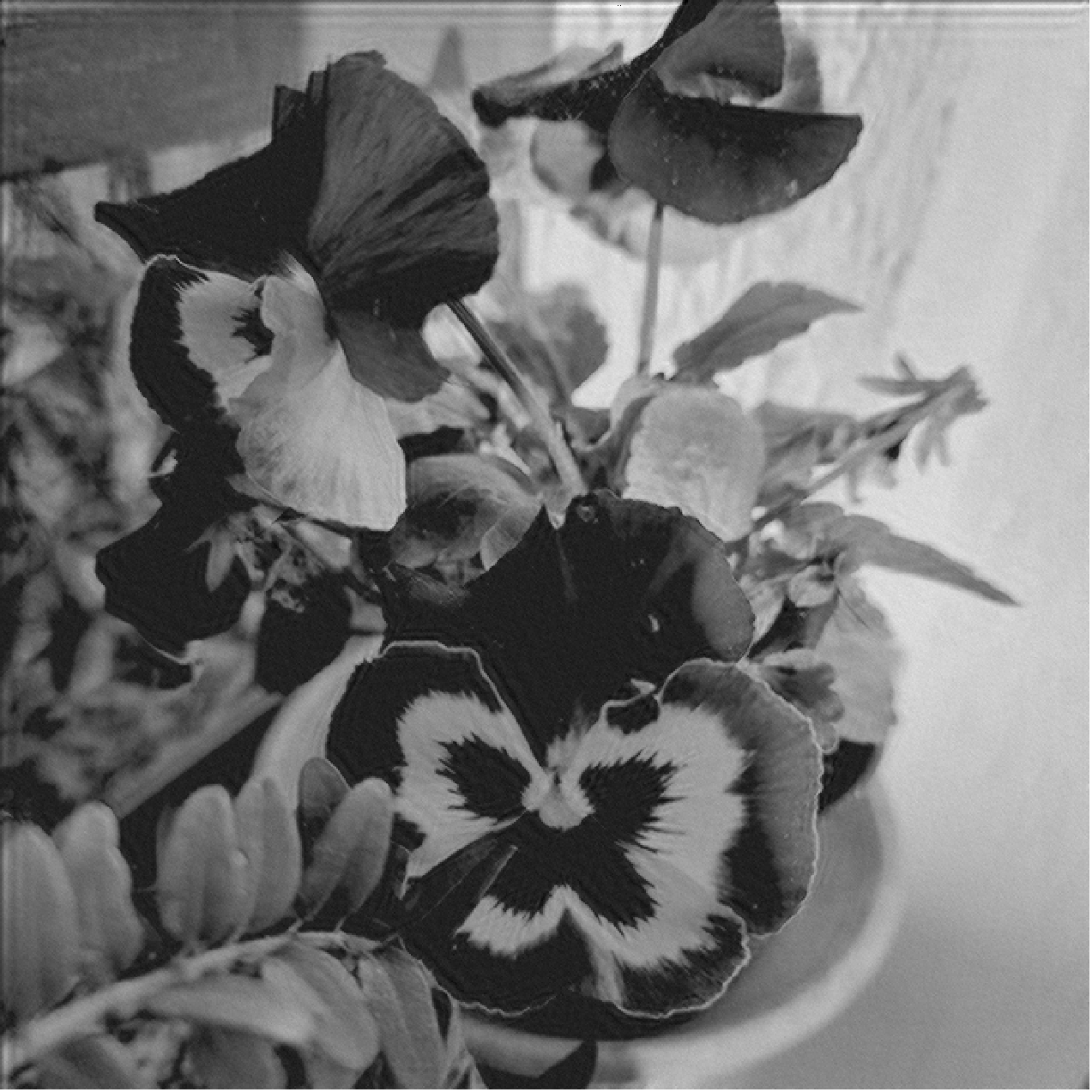}
   \includegraphics[width=0.3\textwidth]{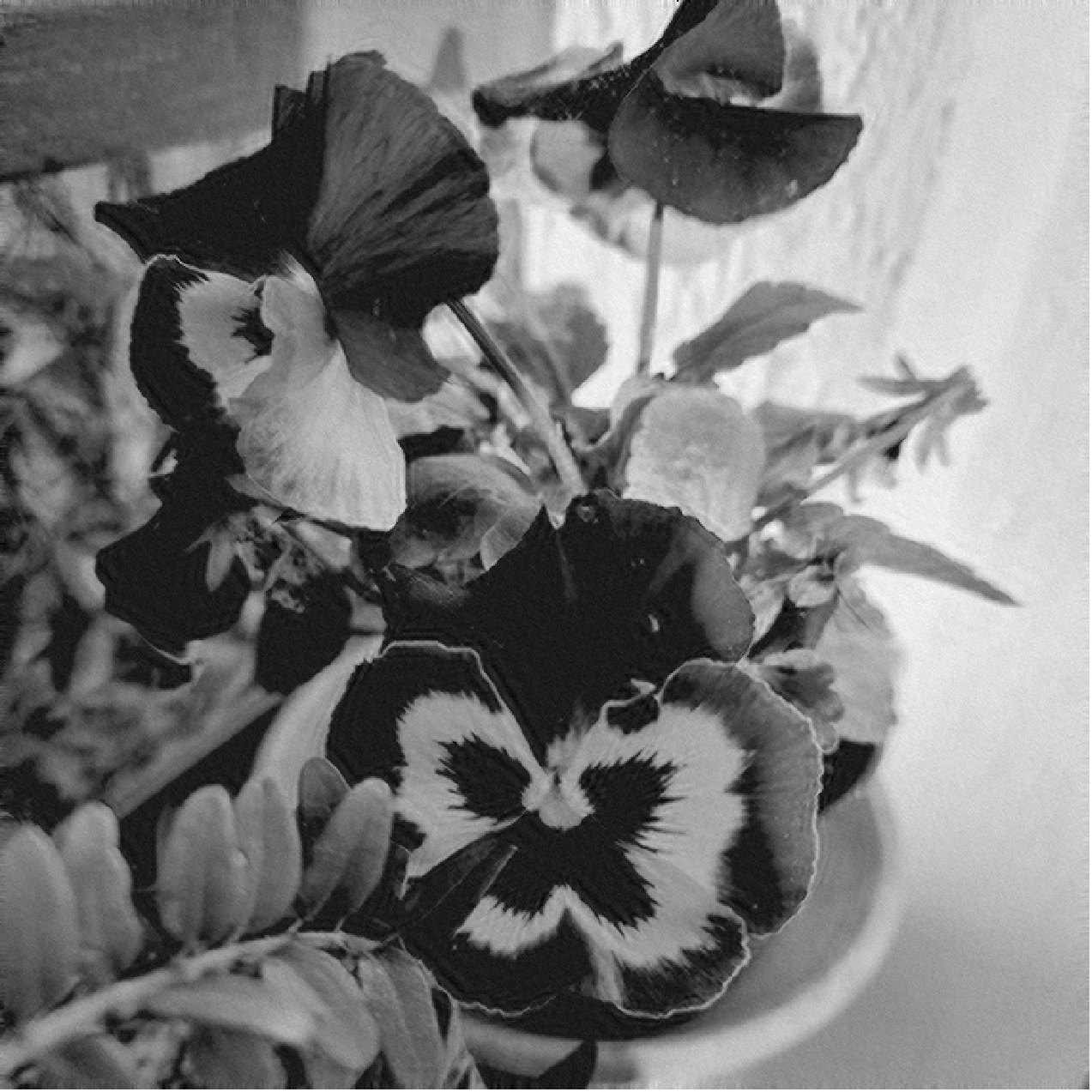}
		\caption{Reconstruction results: true image (top left) compared with reconstructed ones with different BCs and non-flipped/flipped strategies. All the reconstructed images are given by  the discrepancy principle stopping criterion.   Periodic BCs (top middle), zero BCs (top right), reflective BCs (middle left), anti-reflective BCs (middle), periodic BCs with flip (middle right), zero BCs with flip (bottom left), reflective BCs with flip (bottom middle), anti-reflective BCs with flip (bottom right). }
		\label{fig:flowers_images_reconstructed}
	\end{figure}

 Table \ref{tab:results} displays the RRE and PSNR metrics along with their respective iteration counts for two reconstructions: one with the lowest RRE and the other given by the discrepancy principle's stopping rule. Figure \ref{fig:flowers_images_reconstructed} offers a visual comparison of reconstructions achieved via the discrepancy principle under all examined boundary conditions, both with and without the flip operation. The reconstruction using anti-reflective BCs coupled with the flipping strategy yields the most visually appealing result, free from notable artifacts at the edges. The latter impression is further supported by the PSNR values.
 
 {For completeness we conclude the example testing also the MINRES method.  The error curves are shown in Figure \ref{fig:flowers_errors_minres} with a comparison to the GMRES behaviour in the significant cases where the flipped matrices are not symmetric, that is with reflective and anti-reflective BCs. We highlight that the corresponding curves overlap, which is an additional confirmation that the dominant part of the system matrices is the symmetric term as stated in Theorem \ref{th:main}. Since the performances of the two methods are analogous we avoid to include the images reconstructed with the MINRES method.}

\begin{figure}[htbp]
		 \includegraphics[width=0.3\textwidth]{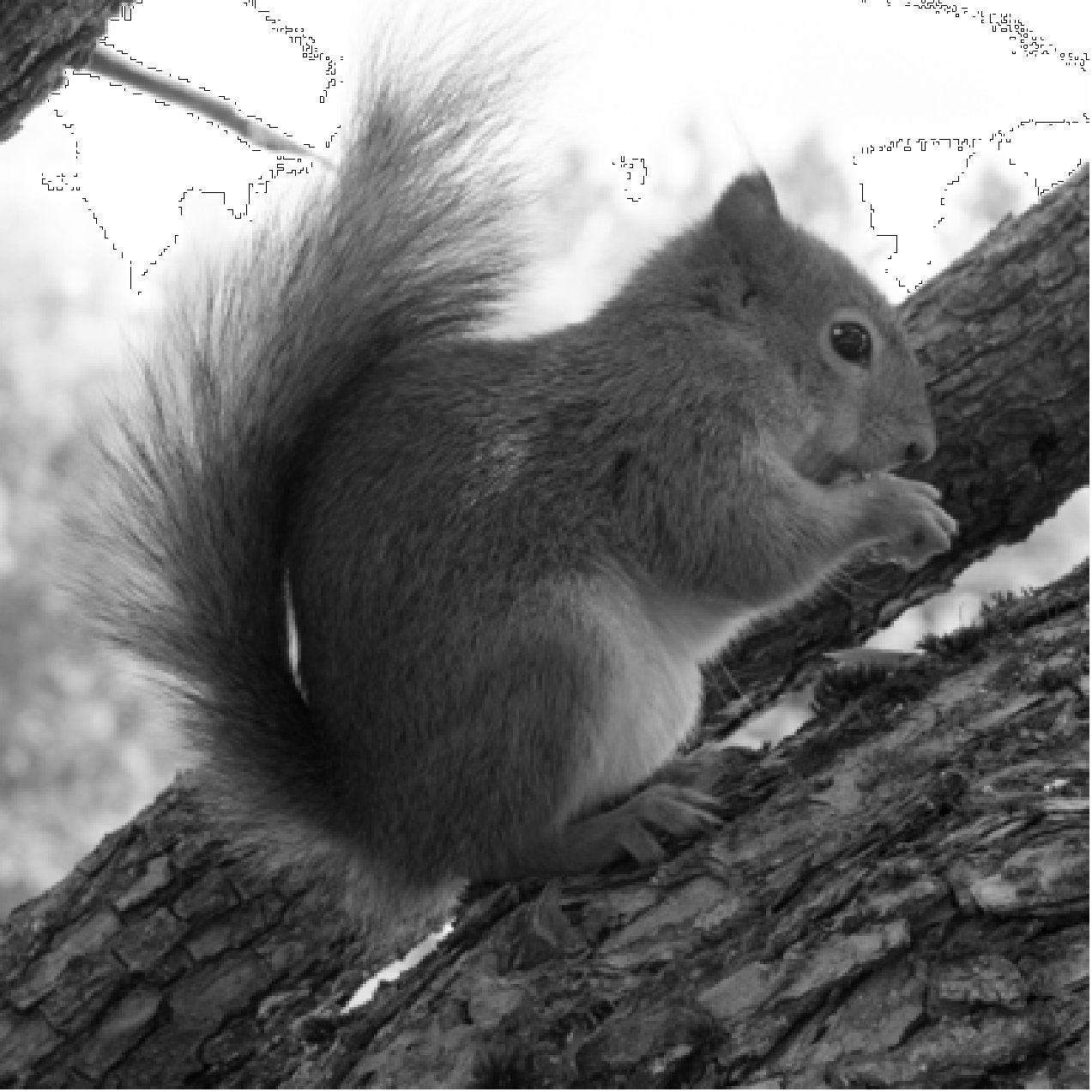}
		 \includegraphics[width=0.3\textwidth]{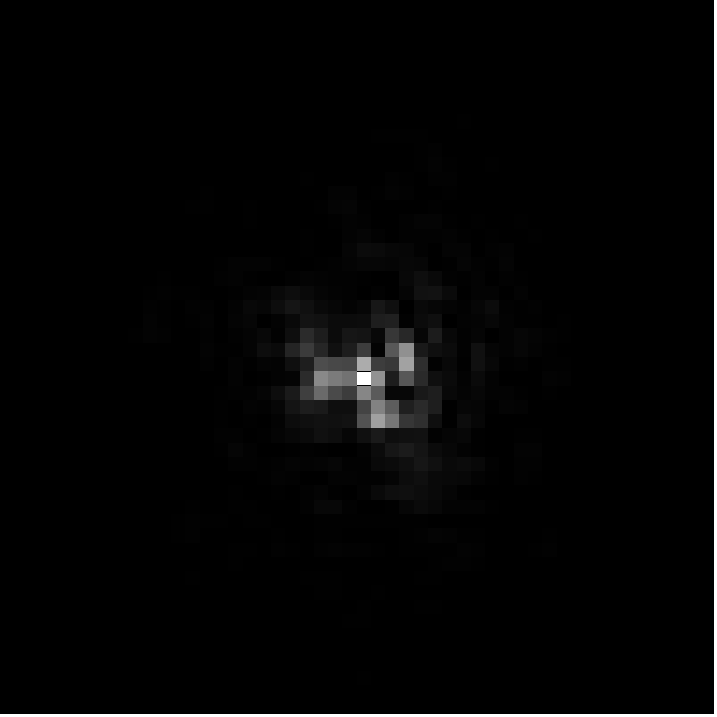}
		 \includegraphics[width=0.3\textwidth]{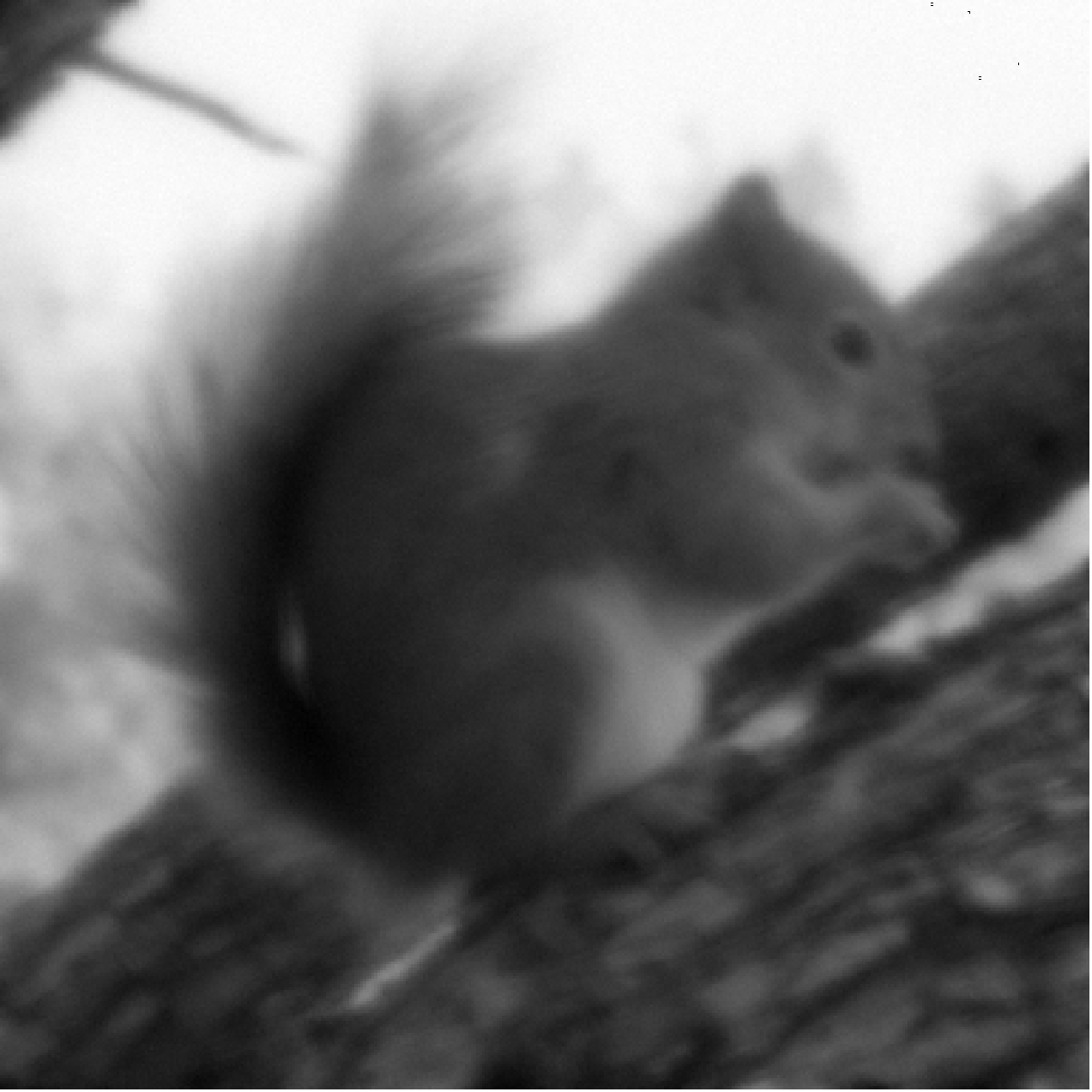}
		\caption{Squirrel image, from left to right: the exact image, the speckle PSF cropped to $[h_{i,j}]_{i,j=-25}^{25}$, and blurred and noisy image.}
		\label{fig:squirrel_images}
	\end{figure}

\paragraph{Squirrel Example}
The second problem that we consider is the reconstruction of an image of size $384\times 384$ pixels representing a squirrel. The image is taken from the TAMPERE17 noise-free image database \cite{TAMPERE17} and resized. In this case, the contamination is by a speckle blur, which arises for example from coherent light interference or atmospheric turbulence, creating a uniformly grainy overlay on the image. The noise level is $\gamma=0.01$. Figure~\ref{fig:squirrel_images} contains the representations of the exact image, the PSF, and the blurred and noisy image.

\begin{figure}[htbp]
		 \includegraphics[width=0.9\textwidth]{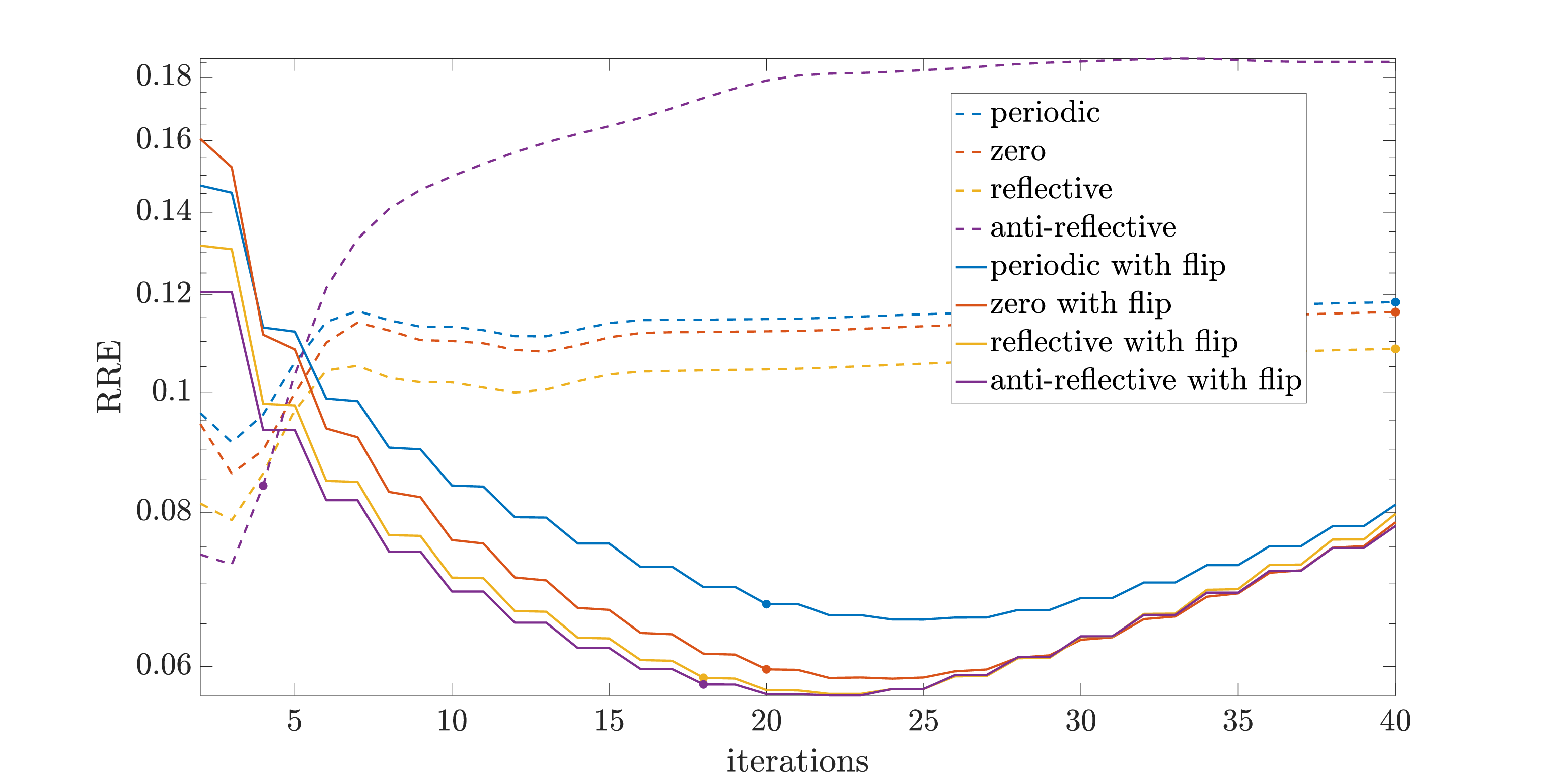}
		\caption{Squirrel image reconstruction: error and iteration behaviours  of the GMRES method applied to the non-symmetrized (dashed lines) and symmetrized (continuous lines) system with different boundary conditions. The dots on the curves mark the iterations where the discrepancy principle is met.}
		\label{fig:squirrel_errors}
	\end{figure}

\begin{figure}[htbp]
		 \includegraphics[width=0.3\textwidth]{squirrel_image_true.eps}
		 \includegraphics[width=0.3\textwidth]{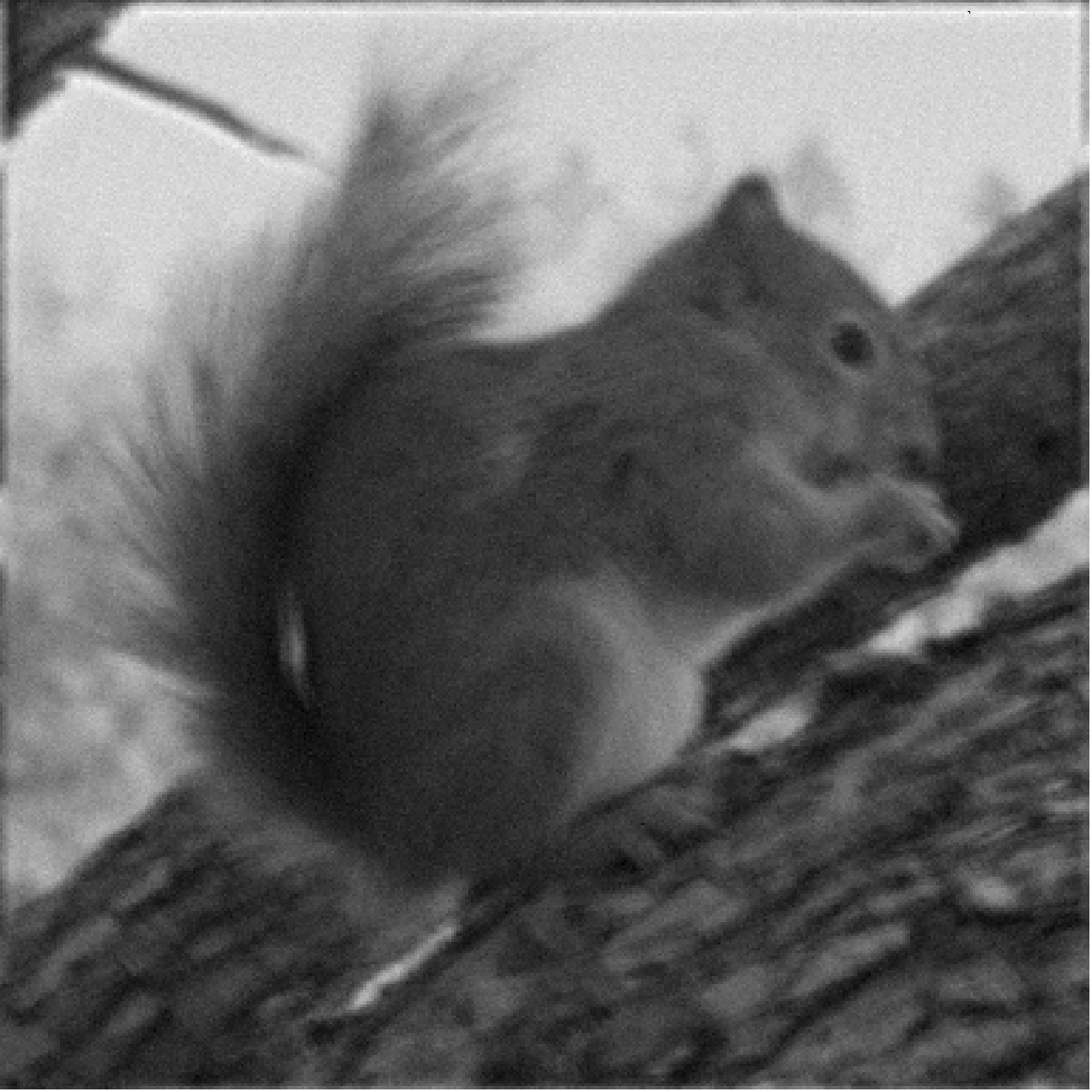}
		 \includegraphics[width=0.3\textwidth]{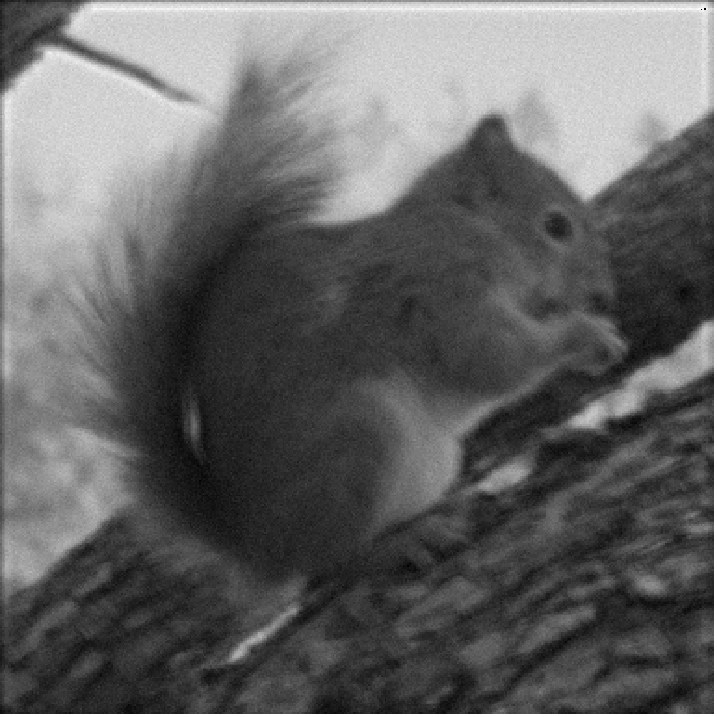}\\
   \includegraphics[width=0.3\textwidth]{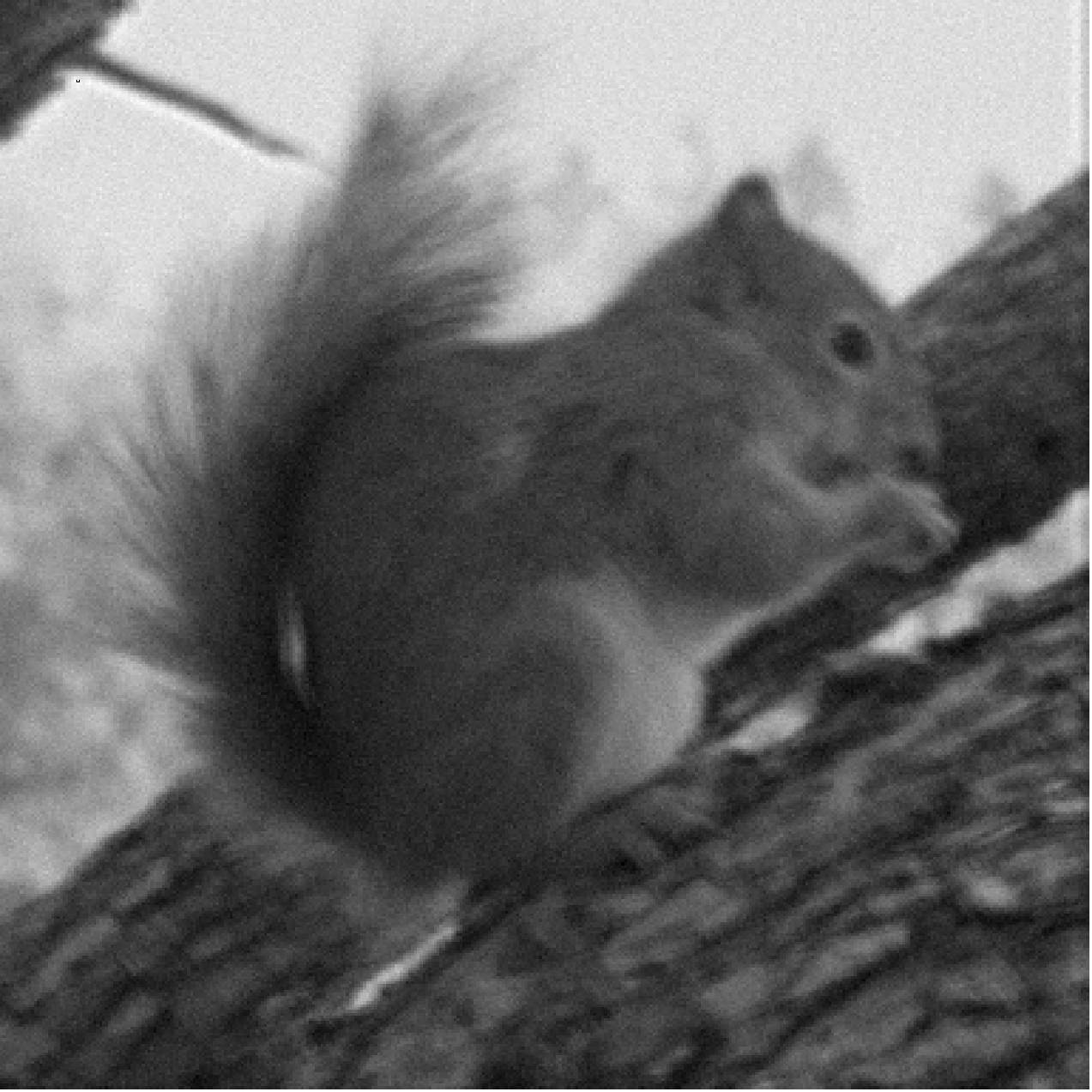}
   \includegraphics[width=0.3\textwidth]{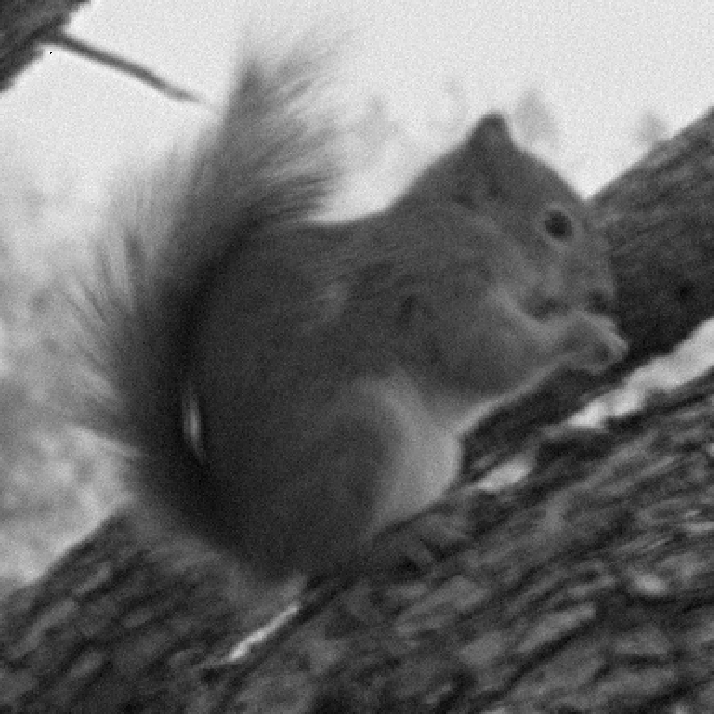}
  		 \includegraphics[width=0.3\textwidth]{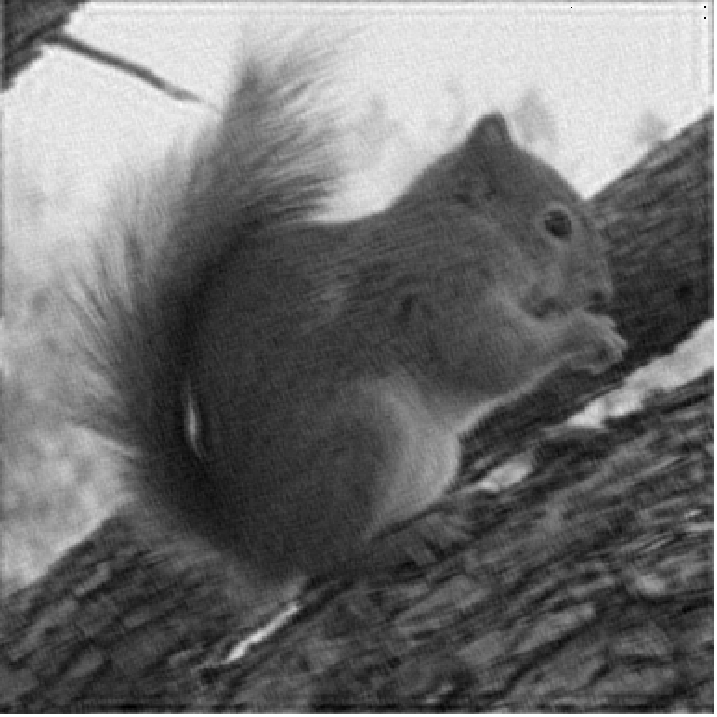}\\
		 \includegraphics[width=0.3\textwidth]{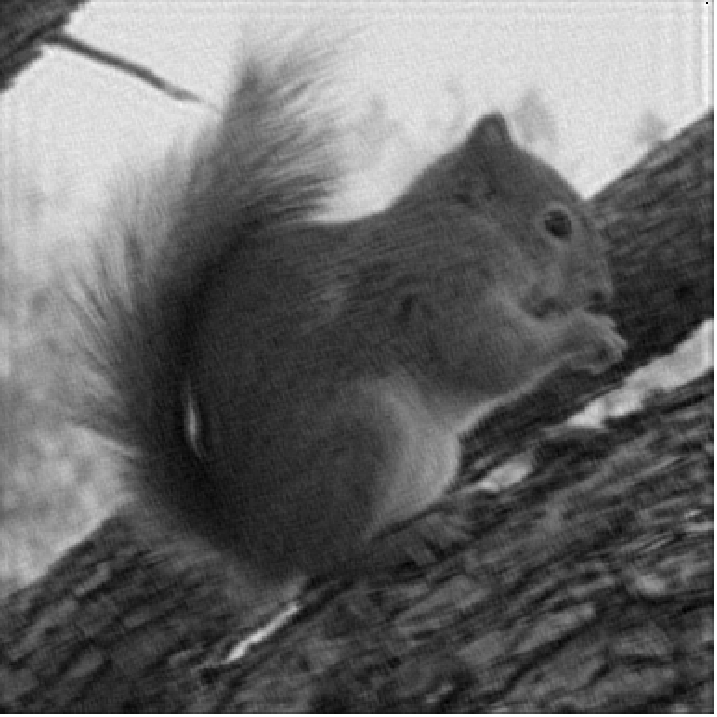}
   \includegraphics[width=0.3\textwidth]{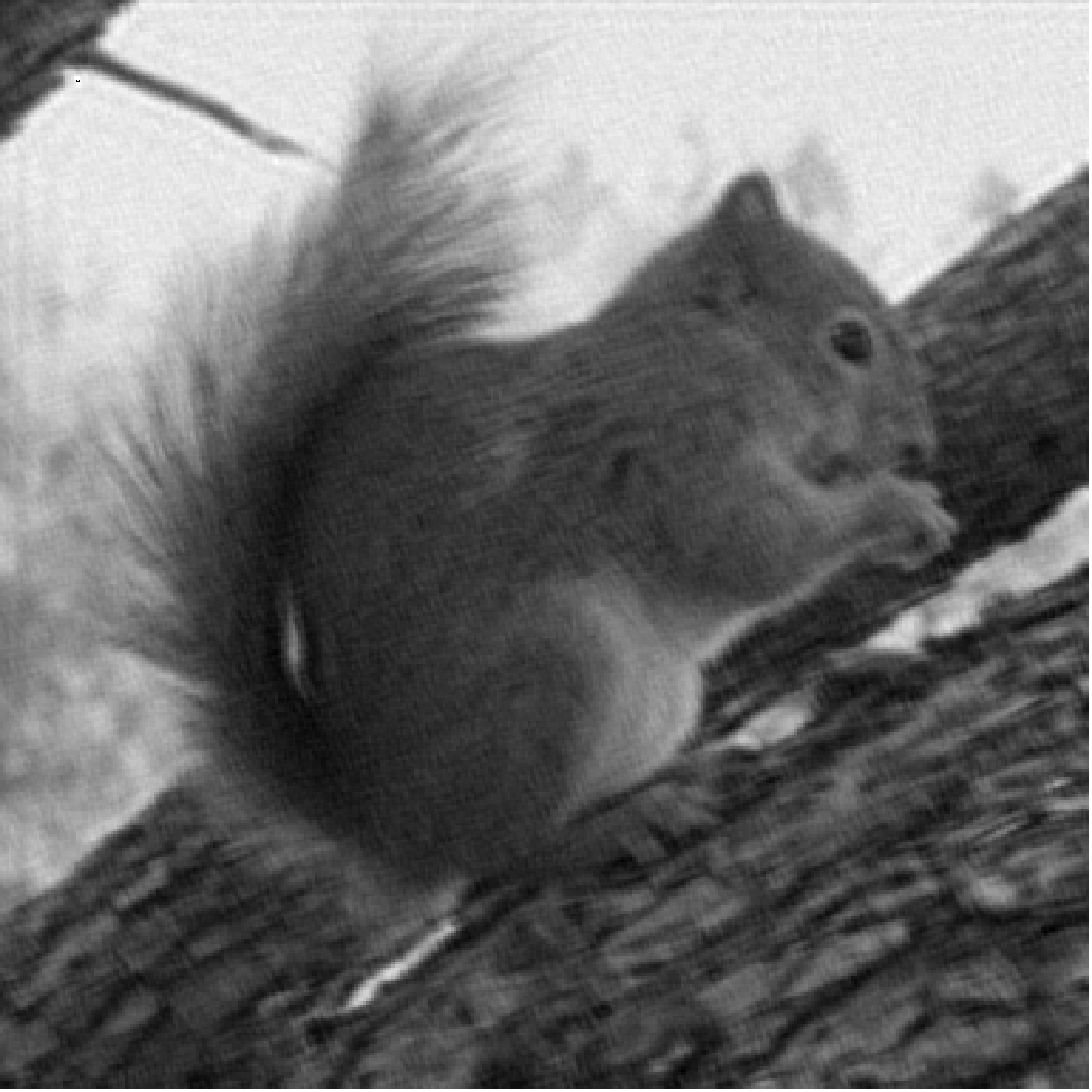}
   \includegraphics[width=0.3\textwidth]{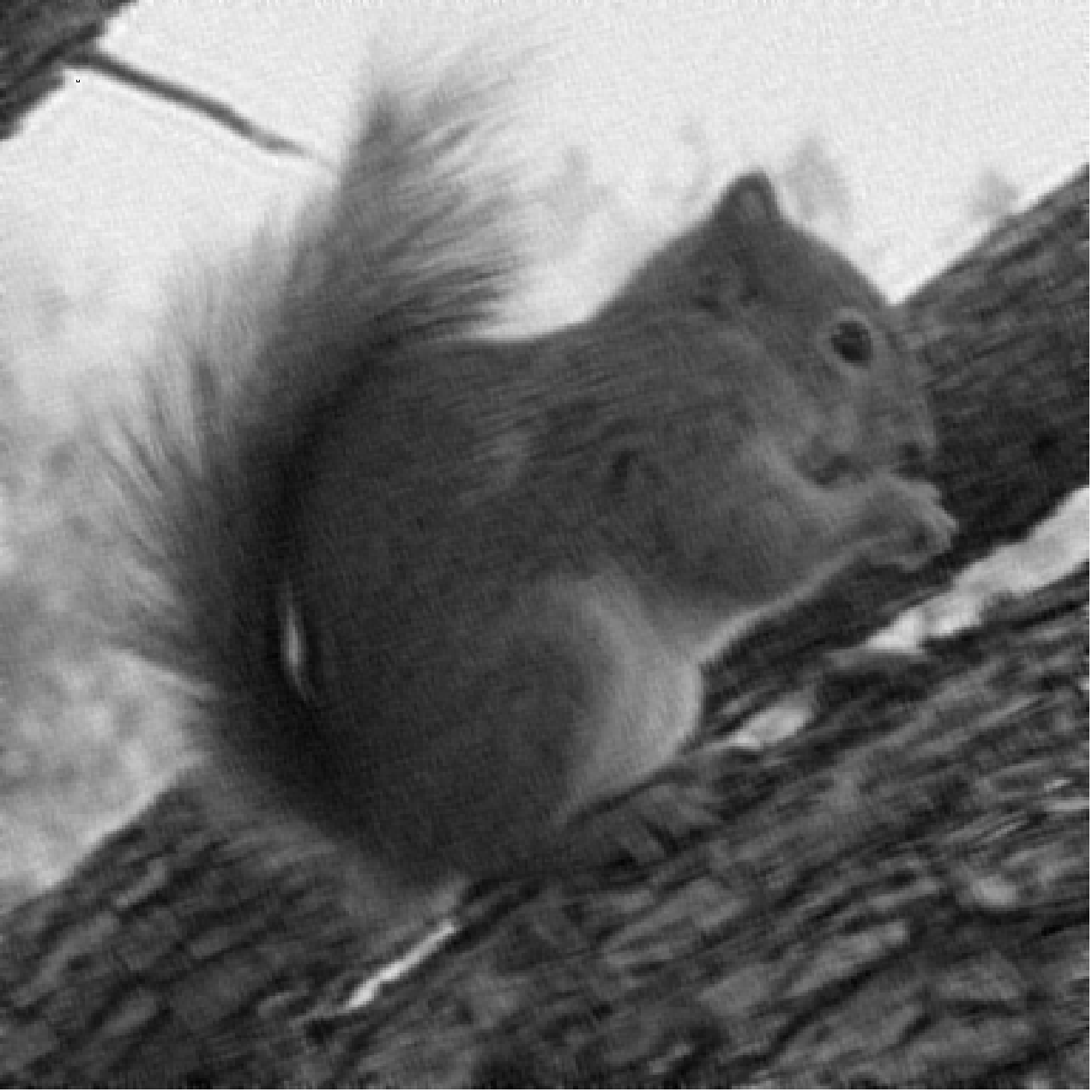}
	\caption{Reconstruction results: true image (top left) compared with reconstructed ones with different BCs and non-flipped/flipped strategies. All the reconstructed images are given by  the discrepancy principle stopping criterion.   Periodic BCs (top middle), zero BCs (top right), reflective BCs (middle left), anti-reflective BCs (middle), periodic BCs with flip (middle right), zero BCs with flip (bottom left), reflective BCs with flip (bottom middle), anti-reflective BCs with flip (bottom right). }
		\label{fig:squirrel_images_reconstructed}
	\end{figure}

 The error behaviour plots in Figure \ref{fig:squirrel_errors} confirm the efficiency of the flipping strategy as a regularising preconditioner. In the bottom part of Table \ref{tab:results} the RRE and PSNR metrics along with their respective iteration counts are shown. Figure \ref{fig:squirrel_images_reconstructed} offers a visual comparison of reconstructions achieved via the discrepancy principle under all boundary conditions, both with and without flipping.
 
 Concerning the selection of BCs, we highlight that in the non-flipped scenario (dashed lines) the discrepancy principle stops the iteration only in the anti-reflective BCs case, which are again the most suitable choice. The second best choice is that of reflective BCs. The advantage of anti-reflective BCs over reflective BCs is less evident after the application of the flipping strategy because the multiplication by $Y_n$ of the operator associated to the anti-reflective BCs results in a ``less symmetric'' operator, as it can be seen by comparing the matrices $W_n^{R}$ and $W_n^{AR}$ in equations \eqref{eq:W_R} and \eqref{eq:W_AR} and {by the analysis in Subsection \ref{ssec:eigenvalue_plots}}.
In summary, when dealing with a significantly non-symmetric PSF and the consequent need of applying the flipping strategy, the decision between reflective and anti-reflective BCs should consider a balance between two key factors: enhancing the reconstruction quality at the boundaries and maintaining the symmetry of the flipped operator. For example, in the previously discussed Flowers example, the use of anti-reflective BCs proved highly effective. Despite the resulting operator being ``less symmetric'' post-flipping, the overall improvement in reconstruction quality was markedly significant. Instead, in the current Squirrel example strategies seem to be equally effective.

 \begin{table}
	\centering
		\begin{tabular}{c|c|c|ccc|ccc}
		\multirow{2}{*}{Example} &\multirow{2}{*}{BCs} &\multirow{2}{*}{Flip}& \multicolumn{3}{c}{\it Best Reconstruction} & \multicolumn{3}{c}{\it Discrepancy Principle}\\
	
	&	& & RRE & PSNR & iter & RRE & PSNR & iter \\
			\hline
		 &Periodic    & no&0.0897 & 26.2495 &  7 & - &-  &- \\
			&Zero   &no &    0.0899 & 26.2248 &  7 & - &-  &- \\
			&Reflective &no   &0.0892 & 26.2929 &  7 & -  &- &- \\   
	\multirow{2}{*}{Flowers}		&Antireflective&no  & 0.0797 & 27.2694 &  7 &-  &-  &-\\
				 &Periodic    & yes& 0.0487 & 31.5548 & 16&    0.0512&  31.1152&  12 \\
			&Zero   & yes& 0.0507 & 31.2087 & 17&    0.0537&  30.7017&  12 \\
			&Reflective &yes   &0.0519 & 31.0003 & 16&    0.0542&  30.6275&  12  \\   
   	&Antireflective&yes  &    {\bf 0.0447} &  {\bf 32.2950} & 17&    {\bf 0.0464}&   {\bf 31.9710} &  12 \\
          \hline
		 &Periodic    & no& 0.0912 & 25.7309  &   3 & - & - & -  \\
			&Zero   &no &    0.0860 & 26.2347  &   3 & - & - & -  \\
			&Reflective &no   &0.0788 & 26.9915  &   3 & - & - & -  \\ 
	\multirow{2}{*}{Squirrel}		&Antireflective& no &     0.0726 & 27.7104  &   3 & 0.0841 & 26.4342 &  4   \\
				 &Periodic    &yes & 0.0655 & 28.6018  & 25 & 0.0674 & 28.3528 & 20  \\
			&Zero   &yes &0.0587 & 29.5598  & 24 & 0.0597 & 29.4086 & 20  \\
			&Reflective & yes  & 0.0570 & 29.8073  & 22 & 0.0588 & 29.5459 & 18  \\
   	&Antireflective&yes  & { \bf 0.0569} & {\bf 29.8309}  & 23 & {\bf 0.0580} &{\bf 29.6535} & 18  \\
		\end{tabular}
		\caption{Values of RRE and PSNR alongside corresponding iteration counts for two different restorations: the one with minimum RRE and the one based on the Discrepancy Principle stopping criterion. The `Flip' column (yes/no) denotes whether the preconditioner $Y_n$ is applied.}\label{tab:results}
\end{table}

\section{Conclusions}\label{sec:end}

Motivated by a recent work on a preconditioned MINRES for flipped linear systems in imaging, we extended the scope of that research for including more precise boundary conditions such as reflective and anti-reflective ones. We proved spectral results for the matrix-sequences associated to the original problem, which justify the use of the MINRES in the current setting. The case of preconditioning has been mentioned in Remark \ref{rem:prec} and it is  essentially left for a future work. The theoretical spectral analysis is supported by a wide variety of numerical experiments, concerning the visualization of the spectra of the original matrices and regarding the convergence speed and regularization features of the considered Krylov methods.

Few open problems remain. Among them, it is worthwhile mentioning 
\begin{itemize}
\item the potential extensions of Theorem \ref{th:main} as mentioned in Remark \ref{rem:extension theorem}, 
\item more precise spectral localization results for the spectra of the original matrix-sequences, 
\item a complete analysis of the resulting matrix-sequences in a preconditioned setting, as suggested in Remark \ref{rem:prec}, also in terms of regularizing features and together with a comparison in terms of computational cost and quality of the signal/image reconstruction,
\item treating the extension for $d=3$ and $m>1$ for specific applications  \cite{qlarger1,qlarger2,qlarger3}: we emphasize that the tools presented in Subsection \ref{ssec:sp-tools} covers completely this setting, even if some adaptations of the results presented in this work are required.
\end{itemize}

\section*{Acknowledgments}
The work of  Paola Ferrari, Isabella Furci and Stefano Serra--Capizzano is partially supported by Gruppo Nazionale per il Calcolo Scientifico (GNCS-INdAM).

Furthermore, the work of Isabella Furci is supported by $\#$NEXTGENERATIONEU (NGEU) and funded by the Ministry of University and Research (MUR), National Recovery and Resilience Plan (NRRP), project MNESYS (PE0000006) – A Multiscale integrated approach to the study of the nervous system in health and disease (DN. 1553 11.10.2022). 

Moreover, 
the work of Stefano Serra-Capizzano was funded from the European High-Performance Computing Joint Undertaking  (JU) under grant agreement No 955701. The JU receives support from the European Union’s Horizon 2020 research and innovation programme and Belgium, France, Germany, Switzerland. 

Finally,
Stefano Serra-Capizzano is grateful for the support of the Laboratory of Theory, Economics and Systems – Department of Computer Science at Athens University of Economics and Business.


\begin{thebibliography}{100}

\bibitem{AR-trasf1}
{\sc A.~Aric\`o, M.~Donatelli, J.~Nagy, and S.~Serra-Capizzano}, {\em The
  anti-reflective transform and regularization by filtering}, in Numerical
  linear algebra in signals, systems and control, vol.~80 of Lect. Notes
  Electr. Eng., Springer, 2011, pp.~1--21.

\bibitem{AR-trasf2}
{\sc A.~Aric\`o, M.~Donatelli, and S.~Serra-Capizzano}, {\em The
  anti-reflective algebra: structural and computational analysis with
  application to image deblurring and denoising}, Calcolo, 45 (2008),
  pp.~149--175.

\bibitem{rearr}
{\sc G.~Barbarino, D.~Bianchi, and C.~Garoni}, {\em Constructive approach to
  the monotone rearrangement of functions}, Expo. Math., 40 (2022),
  pp.~155--175.

\bibitem{BaSe-NLAA}
{\sc G.~Barbarino and S.~Serra-Capizzano}, {\em Non-{H}ermitian perturbations
  of {H}ermitian matrix-sequences and applications to the spectral analysis of
  the numerical approximation of partial differential equations}, Numer. Linear
  Algebra Appl., 27 (2020), p.~e2286.

\bibitem{Book-imaging-astro}
{\sc M.~Bertero and P.~Boccacci}, {\em Introduction to inverse problems in
  imaging}, Institute of Physics Publishing, Bristol, 1998.

\bibitem{Bhatia-book}
{\sc R.~Bhatia}, {\em Matrix Analysis}, vol.~169 of Graduate Texts in
  Mathematics, Springer-Verlag, New York, 1997.

\bibitem{qlarger1}
{\sc J.~F. Cai, R.~H. Chan, and Z.~W. Shen}, {\em A framelet-based image
  inpainting algorithm}, Appl. Comput. Harmon. Anal., 24 (2008), pp.~131--149.

\bibitem{MR1885302}
{\sc D.~Calvetti, B.~Lewis, and L.~Reichel}, {\em On the choice of subspace for
  iterative methods for linear discrete ill-posed problems}, Int. J. Appl.
  Math. Comput. Sci., 11 (2001), pp.~1069--1092.

\bibitem{MR2376196}
{\sc R.~H. Chan and X.-Q. Jin}, {\em An introduction to iterative Toeplitz
  solvers}, vol.~5 of Fundamentals of Algorithms, Society for Industrial and
  Applied Mathematics (SIAM), Philadelphia, PA, 2007.

\bibitem{Chan:1996}
{\sc R.~H. Chan and M.~K. Ng}, {\em Conjugate gradient methods for {T}oeplitz
  systems}, SIAM Rev., 38 (1996), pp.~427--482.

\bibitem{qlarger2}
{\sc V.~Del~Prete, F.~Di~Benedetto, , M.~Donatelli, and S.~Serra-Capizzano},
  {\em Symbol approach in a signal-restoration problem involving block
  {T}oeplitz matrices}, J. Comput. Appl. Math., 272 (2014), pp.~399--416.

\bibitem{donatelli2022symmetrization}
{\sc M.~Donatelli, P.~Ferrari, and S.~Gazzola}, {\em Symmetrization techniques
  in image deblurring}, Electron. Trans. Numer. Anal., 59 (2023), pp.~157--178.

\bibitem{ES-NLAA}
{\sc S.-E. Estr\"om and S.~Serra-Capizzano}, {\em Eigenvalues and eigenvectors
  of banded {T}oeplitz matrices and the related symbols}, Numer. Linear Algebra
  Appl., 25 (2018), p.~e2137.

\bibitem{fasino-tilli}
{\sc D.~Fasino and P.~Tilli}, {\em Spectral clustering properties of block
  multilevel {H}ankel matrices}, Linear Algebra Appl., 306 (2000),
  pp.~155--163.

\bibitem{Ferrari2019}
{\sc P.~Ferrari, I.~Furci, S.~Hon, M.~A. Mursaleen, and S.~Serra-Capizzano},
  {\em The eigenvalue distribution of special 2-by-2 block matrix-sequences
  with applications to the case of symmetrized {T}oeplitz structures}, SIAM J.
  Matrix Anal. Appl., 40 (2019), pp.~1066--1086.

\bibitem{Ferrari2021}
{\sc P.~Ferrari, I.~Furci, and S.~Serra-Capizzano}, {\em Multilevel symmetrized
  {T}oeplitz structures and spectral distribution results for the related
  matrix sequences}, Electron. J. Linear Algebra, 37 (2021), pp.~370--386.

\bibitem{qlarger3}
{\sc P.~J. S.~G. Ferreira}, {\em The stability of a procedure for the recovery
  of lost samples in band-limited signals}, Signal Process., 40 (1994),
  pp.~195--205.

\bibitem{Gander-PinT1}
{\sc M.~Gander and S.~Vandewalle}, {\em Analysis of the parareal time-parallel
  time-integration method}, SIAM J. Sci. Comput., 29 (2007), pp.~556--578.

\bibitem{Gander-PinT2}
{\sc M.~J. Gander}, {\em 50 years of time parallel time integration. multiple
  shooting and time domain decomposition methods}, Contrib. Math. Comput. Sci.,
  9 (2015), pp.~69--113.

\bibitem{GSI}
{\sc C.~Garoni and S.~Serra-Capizzano}, {\em Generalized locally {T}oeplitz
  sequences: theory and applications. Vol. I}, Springer, Cham, 2017.

\bibitem{GSII}
{\sc C.~Garoni and S.~Serra-Capizzano}, {\em Generalized locally {T}oeplitz
  sequences: theory and applications. Vol. II}, Springer, Cham, 2018.

\bibitem{GLT-Eng}
{\sc C.~Garoni, H.~Speleers, S.-E. Estr\"om, A.~Reali, S.~Serra-Capizzano, and
  T.~J.~R. Hughes}, {\em Symbol-based analysis of finite element and
  isogeometric {B}-spline discretizations of eigenvalue problems: exposition
  and review}, Arch. Comput. Methods Eng., 26 (2019), pp.~1639--1690.

\bibitem{ir}
{\sc S.~Gazzola, P.~C. Hansen, and J.~G. Nagy}, {\em Ir tools: a matlab package
  of iterative regularization methods and large-scale test problems}, Numer.
  Algorithms, 81 (2019), pp.~773--811.

\bibitem{MR2287378}
{\sc L.~Golinskii and S.~Serra-Capizzano}, {\em The asymptotic properties of
  the spectrum of nonsymmetrically perturbed {J}acobi matrix sequences}, J.
  Approx. Theory, 144 (2007), pp.~84--102.

\bibitem{MR1413298}
{\sc M.~Hanke}, {\em Conjugate gradient type methods for ill-posed problems},
  vol.~327 of Pitman Research Notes in Mathematics Series, Longman Scientific
  \& Technical, Harlow, 1995.

\bibitem{Book-imaging}
{\sc P.~Hansen, J.~G. Nagy, and D.~P. O'Leary}, {\em Deblurring images}, vol.~3
  of Fundamentals of Algorithms, Society for Industrial and Applied Mathematics
  (SIAM), Philadelphia, PA, 2006.
\newblock Matrices, spectra, and filtering.

\bibitem{evol-flipped3}
{\sc S.~Hon, J.~Dong, and S.~Serra-Capizzano}, {\em A preconditioned {MINRES}
  method for optimal control of wave equations and its asymptotic spectral
  distribution theory}, SIAM J. Matrix Anal. Appl., 44 (2023), pp.~1477--1509.

\bibitem{evol-flipped1}
{\sc S.~Hon, P.~Y. Fung, J.~Dong, and S.~Serra-Capizzano}, {\em A sine
  transform based preconditioned {MINRES} method for all-at-once systems from
  constant and variable-coefficient evolutionary {PDE}s}, Numer. Alg.,  (2023).

\bibitem{Hon2019}
{\sc S.~Hon, M.~A. Mursaleen, and S.~Serra-Capizzano}, {\em A note on the
  spectral distribution of symmetrized {T}oeplitz sequences}, Linear Algebra
  Appl., 579 (2019), pp.~32--50.

\bibitem{evol-flipped2}
{\sc S.~Hon and S.~Serra-Capizzano}, {\em A block {T}oeplitz preconditioner for
  all-at-once systems from linear wave equations}, Electron. Trans. Numer.
  Anal., 58 (2023), pp.~177--195.

\bibitem{MazzaPestana2018}
{\sc M.~Mazza and J.~Pestana}, {\em Spectral properties of flipped {T}oeplitz
  matrices and related preconditioning}, BIT, 59 (2019), pp.~463--482.

\bibitem{MazzaPestana2021}
{\sc M.~Mazza and J.~Pestana}, {\em The asymptotic spectrum of flipped
  multilevel {T}oeplitz matrices and of certain preconditionings}, SIAM J.
  Matrix Anal. Appl., 42 (2021), pp.~1319--1336.

\bibitem{MR2108963}
{\sc M.~K. Ng}, {\em Iterative methods for {T}oeplitz systems}, Numerical
  Mathematics and Scientific Computation, Oxford University Press, New York,
  2004.

\bibitem{MR1718798}
{\sc M.~K. Ng, R.~H. Chan, and W.-C. Tang}, {\em A fast algorithm for
  deblurring models with {N}eumann boundary conditions}, SIAM J. Sci. Comput.,
  21 (1999), pp.~851--866.

\bibitem{MR0383715}
{\sc C.~C. Paige and M.~A. Saunders}, {\em Solutions of sparse indefinite
  systems of linear equations}, SIAM J. Numer. Anal., 12 (1975), pp.~617--629.

\bibitem{MR3323542}
{\sc J.~Pestana and A.~J. Wathen}, {\em A preconditioned {MINRES} method for
  nonsymmetric {T}oeplitz matrices}, SIAM J. Matrix Anal. Appl., 36 (2015),
  pp.~273--288.

\bibitem{TAMPERE17}
{\sc M.~Ponomarenko, N.~Gapon, V.~Voronin, and K.~Egiazarian}, {\em Blind
  estimation of white {G}aussian noise variance in highly textured images},
  Electron. Imaging, 30 (2018), pp.~382--382.

\bibitem{MR1410714}
{\sc S.~Serra}, {\em Preconditioning strategies for {H}ermitian {T}oeplitz
  systems with nondefinite generating functions}, SIAM J. Matrix Anal. Appl.,
  17 (1996), pp.~1007--1019.

\bibitem{taud2}
{\sc S.~Serra-Capizzano}, {\em Spectral behavior of matrix sequences and
  discretized boundary value problems}, Linear Algebra Appl., 337 (2001),
  pp.~37--78.

\bibitem{MR2045058}
{\sc S.~Serra-Capizzano}, {\em A note on antireflective boundary conditions and
  fast deblurring models}, SIAM J. Sci. Comput., 25 (2003), pp.~1307--1325.

\bibitem{ST-ineq}
{\sc S.~Serra-Capizzano and P.~Tilli}, {\em On unitarily invariant norms of
  matrix-valued linear positive operators}, J. Inequal. Appl., 7 (2002),
  pp.~309--330.

\bibitem{MR1366576}
{\sc E.~E. Tyrtyshnikov}, {\em A unifying approach to some old and new theorems
  on distribution and clustering}, Linear Algebra Appl., 232 (1996), pp.~1--43.

\end{thebibliography}
\end{document}